\documentclass{amsart}

\usepackage{amsmath}
\usepackage{amssymb}
\usepackage{amscd}
\usepackage{amsthm}

\numberwithin{equation}{section}
 
\theoremstyle{plain}

\newtheorem{Th}[equation]{Theorem}
\newtheorem{Prop}[equation]{Proposition}
\newtheorem{Le}[equation]{Lemma}
\newtheorem{Cor}[equation]{Corollary}

\theoremstyle{remark}
\newtheorem{remark}[equation]{Remark}

\theoremstyle{definition}

\newtheorem{Def}[equation]{Definition}
\def\e{\emph}
\def\i{\infty}
\def\p{\partial}

\newcommand{\ol}{\overline}
\def\ra{\rightarrow}

\def\Ga{\Gamma}
\newcommand{\E}{\mathbb E}

\begin{document}


\title{Tits Alternative for Closed Real Analytic 4-manifolds of Nonpositive Curvature}

\author{Xiangdong Xie}
\address{Department of Mathematics, Washington University,
St.Louis, MO 63130.}
\email{xxie@math.wustl.edu}

\thanks{2000 \e{Mathematics Subject Classification.} Primary 53C20, 20F67; Secondary 20F65.}

\thanks{\e{Key words and phrases.} {Tits alternative, rank one isometry, Hadamard space,
                      nonpositive curvature,        Tits metric.}}

\vspace{-5mm}



\pagestyle{myheadings}

\markboth{{\upshape Xiangdong Xie}}{{\upshape   Tits Alternative for Analytic 4-manifolds  }}

      
\vspace{3mm}

\begin{abstract}{We study subgroups of fundamental groups of real analytic closed  $4$-manifolds
 with nonpositive sectional curvature. In particular, we are interested in the following question:
 if a subgroup 
of the fundamental group is not virtually free abelian,  does it  contain 
a free  group  of rank two ?  
   The technique involves 
the theory of general metric spaces of nonpositive 
curvature.}
\end{abstract}


\maketitle

\setcounter{section}{0}
\setcounter{subsection}{0}

\section{Introduction}

In this paper we study subgroups of fundamental groups of real analytic closed  $4$-manifolds
 with nonpositive sectional curvature. One would like to know whether  a subgroup   contains  a finite index 
 free abelian   subgroup   if it does not contain a free group of rank two. The same question can be 
 asked in the general setting  of groups acting on Hadamard spaces. We recall 
a Hadamard space is a complete simply connected   metric space with nonpositive curvature
 in the sense of A. D. Alexandrov.   In general there is  the following 
   question (see Section \ref{rank1} for the definition of a proper group action):

\vspace{3mm}

\noindent
\emph{Tits Alternative.} Let $G$ be a group acting   properly and cocompactly
by isometries on a Hadamard space  $X$.   Is it true that every  subgroup $H\subset  G$ 
    contains   either    a  finite index  free abelian  subgroup  or   a free 
group  of rank two?

\vspace{3mm}


The Tits alternative   question  has  been    answered affirmatively      for  the following    spaces: 
    trees (\cite{PV}) or more generally  Gromov hyperbolic  Hadamard spaces  (\cite{G});
    certain cubical complexes (\cite{BSw}); 
   Euclidean buildings of rank $\ge 3$ or symmetric spaces (\cite{T});
     spaces with 
isolated flats (\cite{HR});        
        certain square complexes (\cite{X}).  
The Tits alternative   question  in the general case  appears hard.  In particular it is still open
  for Hadamard $4$-manifolds  and  $CAT(0)$ $2$-complexes.  It   is   not even known (\cite{Sw}) 
 whether $G$ has  an infinite  subgroup where each element is of finite order. 
S. Adams and W. Ballmann (\cite{AB}) showed any amenable subgroup of the group $G$  contains 
  a  finite index  free
 abelian   subgroup.  When  $X$ is  a piecewise smooth $2$-complex where
 each edge is contained in at least two $2$-cells,  W. Ballmann and M. Brin (\cite{BBr}) showed
 either $G$ contains  a free group of  rank two or $X$ is isometric to the Euclidean plane.

Recall  a Hadamard manifold $X$ has higher rank if each geodesic is 
 contained in a $2$-flat,  that is,     a convex subset  of $X$  isometric to 
${\mathbb R}^2$;   $X$ has rank one otherwise. 
When the Hadamard space $X$ is a Hadamard manifold and the group $G$ is torsion free,
       the Tits alternative   can be reduced to  the case when $X$ is 
an irreducible Hadamard manifold, thanks to Eberlein's results (\cite{E1}, \cite{E2})
    on lattices of reducible Hadamard manifolds.   When $X$ is an irreducible Hadamard manifold 
 of higher rank, $X$ is a higher rank symmetric space by the rank rigidity theorem (\cite{B}) and the 
   Tits alternative  
 follows from Tits' theorem.   Thus  one only needs to consider  rank one Hadamard manifolds. 

 The Tits alternative  holds trivially for surfaces. It  is   also not hard to establish it for $3$-manifolds, 
by using some nontrivial results on $3$-manifold topology:

\vspace{3mm}

\noindent
{\bf {Theorem \ref{3d}.}}
\emph{The Tits alternative  holds for subgroups of $\pi_1(M)$, where $M$ is a closed $3$-manifold
 with nonpositive sectional curvature.} 

\vspace{3mm}

   In this paper we mainly consider  the case  of  real analytic closed  $4$-manifolds
 with nonpositive sectional curvature.  Such manifolds have been studied 
by V. Schroeder (\cite{S1}) and C. Hummel and V. Schroeder (\cite{HS1}, \cite{HS2}).
  A higher rank submanifold in a  Hadamard manifold $X$ is a totally geodesic  submanifold which is of
 higher rank as a Hadamard manifold.  When $M=X/\Gamma$ is a real analytic closed  $4$-manifold
 with nonpositive sectional curvature,  one possible class of  higher rank submanifolds  
   have the form $W=Q\times R$, where $Q$ is a nonflat $2$-dimensional  Hadamard manifold. 
 A \emph{cycle} is a finite sequence $W_1, \cdots, W_k$
of distinct higher rank submanifolds of the above form
   such that $W_i\cap W_{i+1}\not=\phi$ ($i=1, 2, \cdots, k-1$) and
 $W_k\cap W_1\not=\phi$.  Below is one of the main results of the paper.

\vspace{3mm}

\noindent
{\bf {Theorem \ref{maint}.}}
\emph{The Tits alternative  holds for subgroups of $\pi_1(M)$  if  $M=X/\Gamma$ is a 
 real analytic closed  $4$-manifold 
 with nonpositive sectional curvature and   there is no cycle in $X$.} 

\vspace{3mm}

There  exist  real analytic closed  $4$-manifolds (see \cite{AS})
satisfying the   assumptions  of   the theorem.




Let $M=X/\Gamma$   be  a  
 real analytic closed  $4$-manifold 
 with nonpositive sectional curvature.   \emph{A singular geodesic }
  in $X$ is a geodesic $c$ of the form  $c=\{q\}\times R\subset Q\times R$, where 
$Q\times R$ is a higher rank submanifold and $Q$ is a nonflat Hadamard $2$-manifold.
 When two higher rank submanifolds  $W_1=Q_1\times R$, $W_2=Q_2\times R$
 intersect,  the intersection  is a $2$-flat of the form   $F=W_1\cap W_2=c_1\times R$ where
$c_1\subset Q_1$ is a geodesic in $Q_1$.    The $R$ directions in $W_1$ and $W_2$ give rise to 
 two parallel families of 
singular geodesics in the $2$ flat $F$. The angle between the two family is a \emph{singular angle}.
 There are only a finite number of singular angles (see Section \ref{real}). 
  U. Abresch   and V. Schroeder (\cite{AS}) have constructed a class of real analytic closed  $4$-manifolds 
 with nonpositive sectional curvature  where the singular angles are all equal to
 $\pi/2$.

\vspace{3mm}

\noindent
{\bf {Theorem \ref{sameangle}.}}
\emph{The Tits alternative  holds for subgroups of $\pi_1(M)$  if  $M=X/\Gamma$ is a 
 real analytic closed  $4$-manifold 
 with nonpositive sectional curvature and    all singular angles are equal to
  $\alpha>\frac{2}{5}\pi$.} 

\vspace{3mm}

Although we are mainly  interested in Hadamard manifolds in this paper, 
  our proof uses general metric spaces with nonpositive curvature  in the sense of A. D. Alexandrov.  
The proof of Theorem \ref{maint} is inspired by the JSJ decomposition in $3$-manifold theory.
Recall a Haken $3$-manifold admits a JSJ decomposition and this decomposition induces a
 graph of groups decomposition for the fundamental group. It follows that the fundamental group
 acts on the Bass-Serre tree  associated to the graph of groups.  Notice that the Bass-Serre tree 
 is not a Hadamard manifold and in general is not locally  compact.  In a similar way we decompose 
  a real analytic  closed $4$-manifold $M$ with nonpositive sectional curvature and construct a 
 $2$-complex associated to the decomposition.  
The fundamental group of $M$ acts on the $2$-complex by isometries. 
The $2$-complex is a $CAT(-1)$ space (see Section \ref{catk} for definition) and should be considered as 
 an analogue of the Bass-Serre
 tree.

The existence of free   subgroups is closely related to the existence of rank one 
isometries.  A  rank one isometry is a hyperbolic isometry (see Section  \ref{hada}   for   definition)  $g$ 
of a Hadamard space $X$ such that 
 no  axis of $g$  bounds a flat half plane, where a flat half plane is a convex  subset 
 of $X$ isometric to the upper half plane: $\{(x,y)\in {\mathbb R}^2: y \ge 0\}$.  
Each isometry of $X$ induces a homeomorphism of the geometric boundary of $X$. 
 If $X$ is a $CAT(-1)$ space or a locally compact Hadamard space,   then 
each  rank one isometry of $X$ acts on the geometric boundary 
  in    the same  way as   a hyperbolic isometry of the real hyperbolic space,
  see Theorem \ref{dyn} or  \cite{B} and \cite{G}.  It follows that   any group  generated by  
two rank one isometries that do  not  share 
 any fixed point in the geometric boundary   contains   a free group of rank two. 
  Notice that  it   is   not necessary  to assume the action is proper.

  For a general hyperbolic  isometry, 
 the  dynamics of the induced homeomorphism on the geometric boundary
have been studied by V. Schroeder (\cite{BGS})  in the Hadamard manifold case 
 and   by K. Ruane (\cite{R}) for locally compact  Hadamard spaces. 
 Let $g$ be a hyperbolic isometry of a locally compact 
   Hadamard space $X$ and $P(g)$ the parallel set of an axis of $g$, that is,
  $P(g)$ is the union of all geodesics that are parallel to an axis of $g$. Denote    by  $\partial_\infty X$ the 
   geometric boundary of $X$.
 $P(g)$ is closed and convex in $X$
 and  $\partial_\infty P(g)$ naturally embeds into  $\partial_\infty X$. 
If $c:R\rightarrow X$ is an axis of $g$, 
 denote     by   $g(+\infty)$  and 
$g(-\infty)$  the points in $\partial_\infty X$ determined by the  rays 
  $c_{|[0, \infty)}$ and $c_{(-\infty, 0]}$  respectively.  V. Schroeder and K. Ruane  showed 
 that for any $\xi\in \partial_\infty X-\partial_\infty P(g)$, the accumulation points of 
 the orbit
 $\{g^i(\xi):i\in {\mathbb Z}\}$   lie in $\partial_\infty P(g)$.
 K. Ruane (\cite{R}) further 
    proved the following result:  if  the Tits distance from $\xi\in \partial_\infty X$ to 
 $g(-\infty)$  is greater than $\pi$,   then there is a neighborhood $U$ of $\xi$
so that $U$ is attracted to $g(+\infty)$ under the iteration of $g$. 
Here we still use  $g$ to denote the homeomorphism 
of the geometric boundary induced by $g$.
  Still the following question remains:

\vspace{3mm}

\noindent
{\bf{Question 1}.}
\emph{With  the above  notation. Is the following   statement  true:  for any compact  subset 
$K\subset  \partial_\infty X-\partial_\infty P(g)$, and  any neighborhood $V$ 
  of $\partial_\infty P(g)$,  there is a positive integer $N$ such that
 $g^n(K)\subset V$ for all $n\ge N$.}  

\vspace{3mm}

If the answer to  Question 1 is yes, then any  group generated by 
two hyperbolic isometries  $f$, $g$  contains  
 a free group of rank two  whenever  $\partial_\infty P(f)$  and $\partial_\infty P(g)$
(considered as subsets of 
  $\partial_\infty X$)   have  empty intersection.

The  existence of rank one isometries in a group is closely related to the  dynamics of the  
group action on the limit set.  
While a lot is known about the group action on the limit set of  
discrete  isometry  groups of  the real hyperbolic space, little is known  
for  groups acting on Hadamard spaces or even Hadamard manifolds. 
  Let $G$ be a  group acting  properly on a Hadamard space such that the limit set 
(see Section \ref{action} for definition)   
 of $G$ contains more than two points. 
   The limit set is closed and  $G$-invariant. If  $X$ is the real hyperbolic space,
   then the limit set is the only nonempty closed and $G$-invariant subset of the limit set.
    This is no longer the case in general, for instance, when 
  $G=G_1\times G_2$ and $X=X_1\times X_2$  where 
$G_1\subset Isom(X_1)$,
  $G_2\subset Isom(X_2)$.  It would be interesting to know when there are more than one 
closed and $G$-invariant subset. 
Recently W. Buyalo and W. Ballmann (\cite{BB})  did some
 interesting work on the topic.  In particular, using their arguments we can show the 
following:

\vspace{3mm}

\noindent
{\bf {Corollary} \ref{diam}.}
\emph{Suppose  $\Gamma$   is   a  group of isometries of a locally compact Hadamard space $X$.  If 
$\Gamma$   does not contain any  rank one isometry,  then $\Lambda(\Gamma)$ 
   has diameter at most $2\pi$ in the Tits   metric.}

\vspace{3mm}


The paper is organized as follows.  In Section \ref{pre} we recall basic facts about Hadamard spaces,
  and collect  results that shall be needed later on.  The topics covered in this section include:
 higher rank submanifolds in the universal covers of closed real analytic $4$-manifolds with nonpositive sectional 
curvature and structure of Tits boundary of such manifolds, rank one isometry and free groups, action on the limit set.
In Section \ref{3manifold} we establish the Tits alternative  for closed $3$-manifolds with nonpositive sectional curvature.
In Section \ref{deco}  we decompose closed real analytic $4$-manifolds, construct the associated 
$2$-complex and use the $2$-complex  to establish 
 the main result of the paper (Theorem \ref{maint}).  In Section \ref{cyclesof} we 
  discuss the Tits alternative       without assuming the   nonexistence  of  cycles. 

\noindent
\textbf{Acknowledgment.} \textit{I am grateful  to Quo-Shin Chi  for  numerous 
  discussions. I would 
  also  like to thank    Bruce Kleiner,   Blake Thornton, Kevin Scannell, 
     Rachel  Roberts, Victor Schroeder, Michael  Kapovich and Alan Reid
 for  remarks and communications on the subject.}

\section{Preliminaries}\label{pre}

In  this section   we recall basic facts about Hadamard spaces,
  and collect  results that shall be needed later on.
We refer the reader to \cite{B}, \cite{BGS}, \cite{BH},  \cite{S1}   and  
  \cite{HS1}    for more details on the material in this section.

\subsection{$CAT(\kappa)$ Spaces} \label {catk}

 Although we are mainly interested in Hadamard manifolds in this paper, general metric spaces 
with upper curvature bounds will play an important  role in the proof of our results. So here we recall the 
basic definitions below.

Let $(X,d)$  be  a 
 metric space.      
A \emph{geodesic}  in $X$   is a continuous map
    $\alpha: I\rightarrow X$   such that, for any point $t\in I$, there exists a neighborhood $U$ of $t$
  with  $d(\alpha(s_1),\alpha(s_2))=|s_1-s_2|$ for   all  $s_1,s_2 \in U$.   
If
 the above equality holds for all $s_1, s_2 \in I$,
then we call $\alpha $ a \emph{minimal geodesic}.
 The image of a geodesic shall also be called a geodesic.   When $I$ is a closed interval
$[a,b]$, we say  $\alpha$ is a geodesic segment of length $b-a$ and  $\alpha$ connects  $\alpha(a)$  
  and  $\alpha(b)$.  
A metric space $X$ is called a \emph{geodesic metric space} if for any  two points  $x,y\in X$ there is a 
   minimal geodesic    segment  connecting them.  

A \emph{ triangle}  in  a metric 
space $(X,d)$  is the union of   three geodesic segments $\alpha_i: [a_i, b_i]\rightarrow X$ ($i=1, 2, 3$)
where  $\alpha_1(b_1)=\alpha_2(a_2)$, $\alpha_2(b_2)=\alpha_3(a_3)$   and  
$\alpha_3(b_3)=\alpha_1(a_1)$.    For any real number $\kappa$,   let $M^2_{\kappa}$ stand for the 2-dimensional simply connected complete Riemannian manifold with constant curvature $\kappa$, and   $D(\kappa)$  denote  the diameter of   $M^2_{\kappa}$ ($D(\kappa)=\infty $ if $\kappa\le 0$).    Given a triangle $\Delta=\alpha_1\cup\alpha_2\cup \alpha_3$
in   a metric space $X$ where $\alpha_i: [a_i, b_i]\rightarrow X$ ($i=1, 2, 3$),    a triangle $\Delta'$ 
in  $M^2_{\kappa}$     is a  \emph{comparison triangle}    for $\Delta$ if they have the same edge lengths, that is,  if  $\Delta'=\alpha_1'\cup\alpha_2'\cup \alpha_3'$
  and  $\alpha_i': [a_i, b_i]\rightarrow  M^2_{\kappa}  $  ($i=1, 2, 3$).   A point $x'\in \Delta'$ corresponds to a point $x\in  \Delta$ if there is some  $i$ and  some $t_i\in [a_i, b_i]$ with $x'=\alpha_i'(t_i)$ and  $x=\alpha_i(t_i)$.
We notice if the perimeter of a triangle $\Delta=\alpha_1\cup\alpha_2\cup \alpha_3$
in $X$ is less than $2D(\kappa)$, that is, if $length(\alpha_1)+length(\alpha_2)+
length(\alpha_3) < 2D(\kappa)$,   then   there is a unique comparison
    triangle (up to isometry) in $M^2_{\kappa}$   for $\Delta$.

\begin{Def}\label{d1}
{A complete metric  space $X$ is called  a \emph {$CAT(\kappa)$ space} if\newline
{(i)}   every two points $x_1, x_2\in X$ with $d(x_1, x_2)< D(\kappa)$ are connected  by a   minimal 
   geodesic segment;\newline
{(ii)}   for any triangle  $\Delta$ in $X$ with perimeter less than $2D(\kappa)$ and any two points $x,y\in \Delta$, the inequality  $d(x,y)\le d(x',y')$ holds,  where  $x'$ and $y'$ are the points on a comparison triangle for $\Delta$ corresponding to $x$ and $y$  respectively. \newline
A complete metric space has \emph{curvature $\le \kappa$ } if each point has a 
$CAT(\kappa)$ neighborhood.}  

\end{Def}

 A simply connected complete Riemannian manifold with sectional curvature 
$\le \kappa$ ($\kappa \le 0$) is $CAT(\kappa)$. 
Simplicial metric trees and more generally  $R$-trees are   $CAT(\kappa)$ for any $\kappa$.  
The following lemma follows from the above definition.

\begin{Le}\label{noshort}
{Let $X$ be a $CAT(1)$ space and $\alpha:  [a,b]\rightarrow X$  a geodesic   with  
     $a<b$  and 
$\alpha(a)=\alpha(b)$.  Then the length of $\alpha$ is at least $2\pi$.}
\end{Le}

\subsection{Hadamard Spaces and Their Ideal Boundaries} \label {hada}

  A $CAT(0)$ space is also called a \emph{Hadamard space}.      Let   $(X,d)$    be  a Hadamard space. 
   Then  the distance function $d: X\times X\rightarrow R$ is convex, there is a unique geodesic segment
 between any two points and $X$ is contractible.
For any $x,y\in X$,  $xy$  denotes the unique
geodesic segment connecting $x$ and $y$.
A \emph{ray} starting from $p\in X$ is a geodesic $c: [0,\infty) \rightarrow X$ with 
$c(0)=p$.   Two rays $c_1$ and $c_2$ are \emph{asymptotic} if  
$d(c_1(t), c_2(t))$  is a bounded function   on the interval  $[0,\infty)$.    The 
\emph{ideal boundary} of $X$ 
 is  the set $\partial X$ of asymptotic  classes of rays in $X$.   For any $p\in X$ 
and any  $\xi\in \partial X$, there is a unique ray (denoted    by  $p\xi$) 
  that starts from $p$ and belongs to
$\xi$.   Thus   for any $p\in X$ we can identify $\partial X$ with the set   of rays 
starting from $p$.     Let  $c$  and $c_i$ ($i=1, 2, \cdots $)
   be rays  starting from $p$;
  we say      $\{c_i\}_{i=1}^\i$ 
 converges to  $c$  if  $c_i$ converges to $c $ uniformly
on compact subsets of $[0, \infty)$.   Similarly  for $x_i\in X$ ($i=1, 2, \cdots$), 
     we say 
 $\{x_i\}_{i=1}^\i$ 
converges to $\xi\in X\cup \partial X$  if  
$px_i$ converges to  $p\xi$ uniformly on compact subsets.
 In this way  we define a topology on  $X\cup \partial X$. 
 It  is    easy to check that this topology is independent of the point $p\in X$.  Both this topology 
 and the induced  topology  on $\partial X$   are  called the \emph{cone topology}. 
  The topology on  $X$ induced by the cone topology coincides with the metric topology on $X$. 
  $\partial X$ together with the cone topology 
  is called the \emph{geometric boundary} of $X$, and denoted   by   $\partial_{\infty} X$.
   We set $\ol X=X\cup \partial_{\infty} X$.  


The Tits metric on the ideal boundary is defined as follows. 
Let $\xi,\eta\in \partial_{\infty} X$. For $p\in X$, let $c_1(t)$ and $c_2(t)$ be the
 rays that start  from $p$ and asymptotic to $\xi$ and $\eta$ respectively.  The 
Tits angle        $\angle_T(\xi, \eta)$    between $\xi$ and $\eta$ is given by:
$$\sin(\frac{\angle_T(\xi, \eta)}{2})=
\lim_{t\rightarrow \infty}\frac{d(c_1(t),c_2(t))}{2t}.$$
This definition is independent of the point $p$.    The \emph{Tits metric} $d_T$ is  the path metric induced by
$\angle_T$.  In particular, if $\xi, \eta\in \partial_{\infty} X$ are in different Tits components, then $d_T(\xi, \eta)=\infty$.    $\p_\i X$ with the  Tits 
metric $d_T$ is   called  the \e{Tits boundary} of $X$ and  denoted    by   $\partial_T X$. 
The Tits topology  and the cone topology are generally quite different.


Below we collect some basic facts concerning  the Tits metric. 
 For more details please see \cite{BGS}  and \cite{BH}.
For any geodesic $c: R\rightarrow X$ in a Hadamard space, we   call  the two points 
in $\partial_{\infty} X$ determined by the two rays $c_{|{[0, +\infty)}}$ and  $c_{|{(-\infty,0]}}$
 the endpoints of $c$, and denote them by $c(+\infty)$ and $c(-\infty)$  respectively.

\begin{Prop}\label{Titsb}
{Let $X$ be a Hadamard space and  $\xi,\eta\in  \partial_T X$. \newline
\emph{(i)} $\partial_T X$  is  a  CAT(1) space;\newline
\emph{(ii)} If $X$ is locally compact and  $d_T(\xi, \eta)>\pi$, then there is a geodesic in 
$X$ with $\xi, \eta$ as endpoints;\newline
\emph{(iii)} If   $X$ is locally compact and $d_T(\xi, \eta)<\infty$,   then there
  is a minimal geodesic  in   $\partial_T X$  from $\xi$ to $\eta$.}
\end{Prop} 

Let  $X$ be a Hadamard space and $g: X\rightarrow X$ an isometry of $X$. 
$g$ is called a \emph{hyperbolic isometry} if it translates a geodesic, that is, if there is a geodesic 
$c: R\rightarrow X$ and a positive number $l$ so that 
$g( c(t))=c(t+l)$  for all $t\in R$;    the geodesic $c$ is called an \emph{axis} of $g$.
All the axes of $g$ are parallel, thus  it  makes sense to denote
    $g(+\infty)=c(+\infty)$,  $g(-\infty)=c(-\infty)$.
Recall two geodesics $c_1, c_2:  R\ra X$ are parallel
 if $d(c_1(t), c_2(t))$ is a bounded function over $R$.

  For a hyperbolic   isometry  $g$, let $Min(g)$ be the union of  all the axes
 of $g$.   The subset  $Min(g)\subset X$  is  closed and convex    in  $X$ and 
     splits  isometrically as $Y\times R$,  where each
 $\{y\}\times R$ ($y\in Y$) is an axis of $g$.  
The geometric boundary $\partial_\infty Min(g)$ naturally embeds into 
$\partial_\infty X$.   

Each isometry $g$ of $X$ induces a homeomorphism of the geometric boundary
$\partial_\infty X$, which   we  still  denote   by   $g$. 

\begin{Th}{\emph{(\cite{R})}}\label{hyperbolicf}
{Let $X$ be a Hadamard space and $g$ a hyperbolic isometry of $X$.
   Then the fixed point set of $g$ on $\partial_\infty X$   is 
  $\partial_\infty Min(g)$.}
\end{Th}

\subsection{Real Analytic 4-manifolds With Nonpositive Sectional Curvature} \label {real}

In this section we recall  some facts concerning the universal covers
of  real analytic closed 4-manifolds with nonpositive sectional curvature. 
The reader is referred to   \cite{S1},  \cite{HS1} and  \cite{AS}  for more details. 

Let $X$ be a Hadamard manifold, i.e., 
 a simply connected complete Riemannian manifold with
nonpositive sectional curvature.   A \emph{$k$-flat}  in   $X$ 
  is a  totally geodesic submanifold isometric to the $k$-dimensional Euclidean
space ${\mathbb E}^k$. We say   $X$ has \emph{higher rank} if each geodesic  in $X$ 
is contained in a $2$-flat, and $X$ has \emph{rank 1} otherwise. 
A complete totally geodesic submanifold of $X$ is  a \emph{higher rank submanifold} if it has
higher rank as a Hadamard manifold.   A \emph{maximal higher rank submanifold}   of  $X$
     is a higher rank submanifold that is   maximal with respect to inclusion.

Let $M$ be a closed Riemannian manifold  with nonpositive sectional curvature and $X$ its universal cover.
Then $M=X/\Gamma$ where $\Gamma$ is the group of deck transformations acting on $X$ as isometries.  
A    complete totally geodesic  submanifold  $W$ of $X$ is \emph{closed}  if  $W/{Stab_{\Gamma}W}$ is compact, where
     $Stab_{\Gamma}W=\{\gamma\in \Gamma:  \gamma(W)=W\}$  is the stabilizer of  $W$ in $\Gamma$.   
   We say $M$ is of rank 1 if $X$ is of rank 1.

\begin{Th}{\emph{(\cite{S1})}}\label{subm}
{Let $M=X/\Gamma$  be a rank 1 closed real analytic 4-manifold of nonpositive sectional curvature. 
Then:\newline
\emph{(i)} each maximal higher rank submanifold of $X$ is isometric to one of the following:
$\E^2$, $\E^3$, $Q\times R$, where $Q$ is a nonflat 2-dimensional Hadamard manifold;\newline
\emph{(ii)} each maximal higher rank submanifold is closed;\newline
\emph{(iii)} there are only a finite number of maximal higher rank submanifolds modulo $\Gamma$;\newline
\emph{(iv)} $W_1\cap W_2\cap W_3=\phi$ for any three distinct  $3$-dimensional 
      higher rank submanifolds $W_1$, $W_2$, $W_3$  of $X$.}
\end{Th}

For $M=X/\Gamma$   as in    Theorem \ref{subm},  we let ${\mathcal {W}}$ be the set of  maximal 
higher rank submanifolds of X that are of the form $Q\times R$ for nonflat 2-dimensional 
Hadamard manifolds $Q$.

\begin{Th}{\emph{(\cite{S1})}}\label{inter}
{Let $M=X/\Gamma$  be a rank 1 closed real analytic 4-manifold of nonpositive sectional curvature,
and $W_1$, $W_2$ maximal  higher rank submanifolds of $X$ with  $W_1\cap W_2\not=\phi$.    Then one of the following 
holds:\newline
\emph{(i)}  $W_1$, $W_2$  are both isometric to $\E^2$. In this case $W_1\cap W_2$ is a point;\newline
\emph{(ii)} $W_1,  W_2\in {\mathcal {W}}$. In this case,   $W_1$ and $W_2$ are perpendicular to each other,
$W_1\cap W_2$ is a $2$-flat and $W_1\cap W_2=c_1\times R\subset Q_1\times R=W_1$ 
for a geodesic $c_1$ of $Q_1$.}
\end{Th}

Notice  for a   complete  totally geodesic submanifold $W$ of $X$, the   geometric  boundary
$\partial_{\infty} W$    naturally  embeds into $\partial_{\infty} X$.

\begin{Th}{\emph{(\cite{HS1})}}\label{Titscom}
{Let $M=X/\Gamma$  be a rank 1 closed real analytic 4-manifold of nonpositive sectional curvature,
  and $\partial_T X$ the Tits boundary of $X$.    Assume ${\mathcal {C}}$ is  a connected component 
of $\partial_T X$.  Then exactly one of the following   statements is true:\newline
\emph{(i)} $\mathcal C$   consists of a single point;\newline
\emph{(ii)} $\mathcal C=\partial_{\infty} F$,  where $F$ is a $2$-flat;\newline  
\emph{(iii)} $\mathcal C=\partial_{\infty} F$,  where $F$ is a $3$-flat;\newline  
\emph{(iv)} $\mathcal C=\bigcup_{W\in {\mathcal W}^*}\partial_{\infty} W$,  where $ {\mathcal W^*}\subset {\mathcal {W}}$ 
  is a subset such that  $\bigcup_{W\in {\mathcal W}^*}W$ is a connected component of 
$\bigcup_{W\in {\mathcal W}}W\subset X$;\newline
\emph{(v)}  $\mathcal C$   is isometric to a closed interval    with  length  $<\pi$.}
\end{Th}

For convenience, we  introduce the following definition.

\begin{Def}\label{cycle}
{Let $M=X/\Gamma$  be a rank 1 closed real analytic 4-manifold of nonpositive sectional curvature.
An  \emph{$n$-cycle} in $X$ is a sequence of $n$  distinct higher rank submanifolds $W_1$, $W_2, \cdots, W_n$
  such that all $W_i\in {\mathcal W}$ and $W_i\cap W_{i+1}\not=\phi$ for all $i$, where
  indices are taken modulo $n$.  A cycle in $X$ is an   $n$-cycle for some $n$.}
\end{Def}

\begin{Prop}{\emph{(\cite{HS1})}}\label{no4cycle}
{Let $M=X/\Gamma$  be a rank 1 closed real analytic 4-manifold of nonpositive sectional curvature.
    If  there   is  an   $n$-cycle in $X$, then $n\ge 5$.}
\end{Prop}

For any metric space $Z$ and any 
   $A\subset Z$,  $\epsilon>0$, the $\epsilon$-neighborhood 
   of $A$ is $N_\epsilon(A)=\{z\in Z: d(z,a)<\epsilon \;\; \text{for some }\; a\in A\}$. 
  Let $X$ be a Hadamard manifold.
 Two  complete totally geodesic submanifolds $H_1$  and $H_2$ of  $X$   
 are parallel if there  is some $\epsilon>0$   with  $H_2\subset N_\epsilon(H_1)$
  and $H_1\subset N_\epsilon(H_2)$.
  Let $H$ be a  complete totally geodesic submanifold of $X$,
the parallel set $P_H$ of $H$ is the union of all complete totally geodesic submanifolds of $X$
that are parallel to $H$.  $P_H$ is a closed convex subset of $X$ and splits isometrically as 
 $P_H=H\times Y$, where   each  $H\times \{y\}$ ($y\in Y$)    is   a   
   complete totally geodesic submanifold  
parallel to $H$.   In general $P_H$ has boundary and is not a manifold. 
 When the Riemannian manifold $X$ is real analytic, $P_H$ is  
complete totally geodesic.

\subsection{Rank  One Isometries  and Free Groups } \label{rank1}

Let  $X$ be a Hadamard space. 
      A \emph{flat half  plane}  in $X$  is 
   the  image of an isometric embedding  
  $f:\{(x,y)\in \E^2: y\ge 0\}\rightarrow X$, and in this case we say the geodesic 
  $c:R\rightarrow X$,  
$c(t)=f(t,0)$   bounds the 
flat half plane.   

\begin{Def}\label{rank1iso}
{A hyperbolic isometry $g$ of  a Hadamard space $X$ is  called  a \emph{rank  one isometry}   
if   no  axis of $g$  bounds  a  flat half plane.}
\end{Def}

If $X$ is a   Gromov hyperbolic  Hadamard space,   
    then    each hyperbolic isometry
of $X$ is  rank one.  We note a  $CAT(-1)$ space  is  Gromov hyperbolic. 
 The following theorem is due to W. Ballmann in the case of
 locally compact  Hadamard spaces and to M. Gromov in the case of Gromov hyperbolic spaces.
    Recall an  isometry   $g$ of a Hadamard space  $X$ induces a homeomorphism of   $\ol  X$,
   which we still denote  by  $g$.

\begin{Th}{\emph{(\cite{B}}, \emph{\cite{G})}}\label{dyn}
{Let $X$ be a  Hadamard space that is either locally compact  or  
Gromov hyperbolic,  
   and 
  $g$      a rank one isometry of  $X$. 
   Given any neighborhoods
$U$ of $g(+\infty)$ and $V$ of $g(-\infty)$ in $\ol  X$,  there is an $n\ge 0$ 
such that  $g^k(\ol  X- V)\subset U$ and
$g^{-k}(\ol  X- U)\subset V$  whenever
$k\ge n$.}
\end{Th}

Theorem \ref{dyn} in particular implies $g(+\infty)$ and $g(-\infty)$ are the only fixed points 
of $g$ in $\ol  X$   under the conditions of the theorem.

Let $G$ be a group acting by isometries on a Hadamard space $X$.  
The action is said to be \emph{proper}
if for   any compact subset $K\subset X$   the set 
 $\{g\in G: g(K)\cap K\not=\phi \}$ is finite.
   $G$ also 
acts 
on  $\ol X $     as homeomorphisms.   A subset $A\subset \ol X $ 
is \emph{$G$-invariant} if   $g(A)=A$ for all   $g\in  G$.   
 Theorem  \ref{dyn}   has  the following two corollaries.

\begin{Cor}\label{3possibility}
{Let $X$ be a  Hadamard space that is either locally compact  or  
Gromov hyperbolic,
  $G$     a  group of isometries of  $X$    and 
   $g\in G$    a rank one isometry   with fixed points $g(+\infty)$,  $g(-\infty)$.
 Then one of the following holds:\newline
\emph{(i)} $g(+\infty)$   or   $g(-\infty)$  is fixed by all elements of $G$; \newline
\emph{(ii)}   some axis $c$ 
  of $g$  is  $G$-invariant; \newline 
\emph{(iii)}   $G$ contains  a free group of rank two.}
\end{Cor}






Recall a group is \emph{virtually free abelian} if a finite index subgroup is free abelian.
A   virtually    infinite cyclic   group is similarly defined. 

\begin{Cor}\label{1possibility}
{Let $G $   act    properly   and cocompactly by isometries on a $CAT(0)$ space,  and $H$ 
a  subgroup of $G$.    If $H$ contains a rank one isometry, then $H$ either
  is   virtually    infinite cyclic or contains a free group of rank two.}
\end{Cor}




For any two isometries $g$ and $h$ of a Hadamard space $X$, $<g,h>$ denotes the group generated
by $g$ and $h$.

\begin{Th}{\emph{(\cite{R})}}\label{kim}
{ Let $X$ be a   locally   compact Hadamard space and $g$, $h$ two hyperbolic isometries of $X$. 
If $d_T(\xi, \eta)>\pi$    whenever 
   $\xi\in \{g(+\infty),g(-\infty)\}$ and $\eta\in \{h(+\infty),h(-\infty)\}$, 
   then   $<g,h>$ 
   contains a free   group  of rank two.}
\end{Th}

It follows from  Theorem \ref{kim}  that   if    there are 
two distinct Tits components
$C_1$ and $C_2$ such that 
 $\{g(+\infty), g(-\infty)\}\subset C_1$
   and $\{h(+\infty), h(-\infty)\}\subset C_2$, 
 then  the group 
  $<g,h>$   
contains a free group of rank two.

\subsection{Action on the Limit Set} \label{action}

Let $Isom(X)$ be the group  of   all isometries of a Hadamard space $X$,
   and  $\Gamma\subset Isom(X)$   any subgroup.    A point 
$\xi\in \partial_{\infty} X$ is a \emph{limit point}  of $\Gamma$ if 
there is a sequence of elements  $\{\gamma_i\}_{i=1}^{\infty}\subset \Gamma$
 with $\gamma_i(x) \rightarrow \xi$ for some (hence any) $x\in X$. 
 The \emph{limit set} $\Lambda(\Gamma)\subset \partial_{\infty} X $ 
of $\Gamma$ is the set of limit points of  $\Gamma$. 
  It   is   easy to check that  
$\Lambda(\Gamma)$ is closed (in the cone topology) and $\Gamma$-invariant.

\begin{Def}\label{minimal}
{A     nonempty closed and $\Gamma$-invariant subset $M\subset \partial_{\infty}X$ is  \emph{$\Gamma$-minimal }
if it does not contain any proper subset that is closed and $\Gamma$-invariant.}
\end{Def}

When the Hadamard space $X$ is locally compact, the geometric boundary $\partial_{\infty}X$
    and all its closed subsets   are  compact.  It follows from Zorn's lemma that 
   $\Gamma$-minimal set always exists when $X$ is locally compact.

We recall   two point $\xi, \eta\in \partial_{\infty}X$  are \emph{$\Gamma$-dual}
  if there is a sequence of elements $\{\gamma_i\}_{i=1}^{\infty}\subset \Gamma$
  such that $\gamma_i(x)\rightarrow  \xi$ and  $\gamma_i^{-1}(x)\rightarrow  \eta$  for some (hence any) $x\in X$ 
  as $i\rightarrow \infty$.      Clearly if  $\xi, \eta\in \partial_{\infty}X$  are \emph{$\Gamma$-dual}
then $\xi, \eta\in \Lambda(\Gamma)$.    For any $\xi\in \Lambda(\Gamma)$,    $D_{\xi}\subset \Lambda(\Gamma)$ 
   denotes the set of points that are   $\Gamma$-dual to $\xi$.  It   is  not hard to check that
 $D_{\xi}$ is closed and $\Gamma$-invariant.  

\begin{Le}{\emph{(\cite{B})}}\label{getrank1}
{ Let $X$ be  a locally compact Hadamard space and $\Gamma\subset Isom(X)$ a subgroup.    If 
$\xi, \eta\in \partial_{\infty}X$  are $\Gamma$-dual and $d_T(\xi,\eta)>\pi$, then
   $\Gamma$ contains rank one isometries.}
\end{Le}

\begin{Le}{\emph{(\cite{CE})}}\label{dual}
{Let $X$ be  a Hadamard space and $\Gamma\subset Isom(X)$ a subgroup. 
If $\xi\in \Lambda(\Gamma)$  and $\eta\in \partial_{\infty}X$    are the endpoints of a geodesic in $X$, then 
$D_{\xi}\subset \overline{\Gamma(\eta)}$.}
\end{Le}

The argument in the following proof  belongs to  Ballmann and Buyalo (\cite{BB}), 
    who   assumed     $\Lambda(\Gamma)= \partial_{\infty} X $.

\begin{Prop}\label{dis}
{Suppose  $\Gamma$   is   a  group of isometries of a  locally compact Hadamard space $X$.  If 
$\Gamma$   does not contain any  rank one isometry,   then
$d_T(m, \eta)\le \pi$ for any $\Gamma$-minimal  set $M\subset \Lambda(\Gamma)$   and any 
$m\in M$, $\xi\in \Lambda(\Gamma)$.}
\end{Prop}

\begin{proof}
Let $M \subset \Lambda(\Gamma)$ be an arbitrary $\Gamma$-minimal  set and $m\in M$.
Clearly $M=\overline{\Gamma(m)}$.
Assume there is some $\eta \in  \Lambda(\Gamma)$ with $d_T(m, \eta)>\pi$.
By  Proposition \ref{Titsb}(ii),  there is a geodesic in $X$ with $m$ and $\eta$ as endpoints. 
   Lemma \ref{dual} implies  $D_{\eta}\subset \overline{\Gamma(m)}=M$.  Since
   $D_{\eta}$ is   nonempty,  closed and $\Gamma$-invariant and $M$ is $\Gamma$-minimal ,  $D_{\eta}=M$.
In particular,  $m\in M=D_{\eta}$    and  $m$, $\eta$ are $\Gamma$-dual. 
Now  Lemma \ref{getrank1} implies there are rank one isometries in $\Gamma$
  since   we assumed $d_T(m, \eta)>\pi$.  

\end{proof}

\begin{Cor}\label{diam}
{Suppose  $\Gamma$   is   a  group of isometries of a locally compact Hadamard space $X$.  If 
$\Gamma$   does not contain any  rank one isometry,  then $\Lambda(\Gamma)$ 
   has diameter at most $2\pi$ in the Tits   metric.}
\end{Cor}




\section{ 3-manifold Groups}\label{3manifold}

In this section we   study  $3$-manifold groups and establish the Tits alternative for fundamental groups 
 of closed $3$-manifolds with nonpositive sectional curvature (Theorem \ref{3d}).
 Although Theorem \ref{haken}
should be known  to many people, we still include a sketch of the proof. The main reason for doing so 
     is that
 our proof uses group actions on trees and is  a simplified version of the proof in Section \ref{deco}.
 The reader  is referred to \cite{H} and \cite{K} for   definitions    and  basic facts
      concerning $3$-manifolds.


\begin{Def}\label{ta}
{A group $G$ has the \emph{TA-property} if it is virtually
free abelian  or 
contains a free group of rank two.}
\end{Def}

The following lemma is clear.

\begin{Le}\label{finiteid}
{A group has the TA-property if   some  finite index subgroup  does.} 
\end{Le}

\begin{Th}\label{haken}
{Let $M$ be a  Haken $3$-manifold  whose boundary has zero Euler
 characteristic, and $H\subset \pi_1(M)$  a subgroup. Then $H$ has the TA-property
 except in the following two cases:\newline
\emph{(i)} $M$ is finitely covered by a torus bundle over       $\mathbb{S}^1$ 
   and $H$ has finite index in    $\pi_1(M)$;\newline
\emph{(ii)} $M$ is finitely covered by a   $\mathbb{S}^1$-bundle
      over the torus  and $H$ has finite index in    $\pi_1(M)$.}
\end{Th}

\noindent
{\e{Sketch of proof.}}
By Lemma \ref{finiteid} we  may assume $M$ is orientable by considering its orientable double cover if necessary.
  $M$ admits the so-called JSJ decomposition: there is a collection of disjoint   embedded tori
    $\{T_1, \cdots, T_k\}$   in $M$ such that
    the homomorphism $\pi_1(T_i)\rightarrow \pi_1(M)$ 
 induced by inclusion is injective for each $i$, and each component
 of $M-\cup_{i=1}^k T_i$ is either a Seifert manifold or an  atoroidal manifold.  Let $\pi:\tilde M\rightarrow M$ 
 be the projection from the universal cover to $M$. Each component of
$\cup \pi^{-1}(T_i)$ is homeomorphic to $\mathbb{R}^2$, and is called a \emph{plane}. 
Now we construct a graph $T$ from the induced decomposition  of $\tilde M$.  
   The vertex  set of $T$ is in one-to-one correspondence with the set of components of 
$\tilde M-\cup_{i=1}^k \pi^{-1}(T_i)$.  Two vertices are joined by an edge if the intersection of 
     the closures  of the 
  corresponding 
 components is a plane. The fact that each plane is separating implies the graph $T$ is actually a tree.  
The action of  $\pi_1(M)$ on $\tilde M$     preserves the decomposition and induces an action on the tree $T$.
Hence  any subgroup $H\subset \pi_1(M)$ also acts on  $T$.  

For the action of a group $H$ on a tree $T$: (see \cite{PV}) if $H$ does not
  contain a free group of rank two,  then one of the following holds:
(1) $H$ fixes a point in $T$, 
   (2) $H$ fixes   some    $\xi\in \p_\i T$, 
   (3) $H$ leaves invariant   a geodesic $c$ in  $T$.
  Now  we need to analyze these 
    three cases.

 First suppose  $H$ fixes a point in $T$.  We may assume $H$ fixes a vertex of   $T$ by passing to 
an index two subgroup if necessary.  The component of $\tilde M-\cup_{i=1}^k \pi^{-1}(T_i)$
corresponding to this vertex  projects down   either to a Seifert manifold  $N$ or to 
an atoroidal manifold $N$, and $H$ is a subgroup of $\pi_1(N)$.
   If  $N$ is an   atoroidal manifold, 
Thurston's theorem says $N$ admits a hyperbolic 
 structure with finite volume. In this case it is clear that $H$ 
has the TA-property. 
    If  $N$  is  a  Seifert manifold,  then  $H$ has the TA-property
  unless  $H$ has finite index in    $\pi_1(M)$   and
$M$ is finitely covered by a circle  bundle over the torus.



Suppose   $H$ stabilizes a geodesic $c$ in $T$.   
  After passing to an index two subgroup if necessary, each  $h\in H$ translates  $c$.
  By considering the restricted action of $H$ on $c$ we have an exact sequence:
$1\rightarrow H'\rightarrow H\rightarrow  {\mathbb Z}\rightarrow 1$, where $H'\subset H$ is the subgroup
  consisting of elements of $H$ that fix $c$ pointwise.  Since the  stabilizer  of each edge 
   in  $\pi_1(M)$  is 
       ${\mathbb Z}^2$, $H'$ is 
${\mathbb Z}^k$ with $k\le 2$.  
If $k\le 1$, then $H$ is clearly virtually ${\mathbb Z}^{k+1}$.
If $k=2$,  by consideration of   cohomological dimension we see
$H$ has finite index in    $\pi_1(M)$   and
$M$ is finitely covered by a  torus  bundle over   a circle. 
The case when $H$ fixes some  $\xi\in \p_\i T$ can be handled similarly. 

\qed

\begin{remark}\label{sol}
\emph{In  the  two exceptional cases of   Theorem \ref{haken},   both   $\pi_1(M)$ and $H$ are virtually 
solvable.}
\end{remark}





Our next goal is to establish the    Tits alternative  for fundamental groups 
 of closed $3$-manifolds with nonpositive sectional curvature. 
The following is a special case of E. Swenson's theorem.
  
\begin{Th}{\emph{(\cite{Sw})}}\label{swenson}
{Let $X$ be a Hadamard space and $G$ a group acting  properly    and cocompactly by
 isometries on $X$.  Suppose $H, K\subset G$ are  subgroups  and $A,B\subset X$ are 
closed convex subsets such that $h(A)=A$, $k(B)=B$ for all $h\in H, k\in K$ and $A/H$, 
$B/K$    are  compact.
If $A\cap B\not=\phi$, then $H\cap K$ acts cocompactly on $A\cap B$.}

\end{Th}

  Recall for any group $G$ and any $g\in G$, 
the centralizer  of $g$ in $G$ is   
 $C_g(G)=\{\gamma\in G: \gamma g=g\gamma\}$.

\begin{Th}{\emph{(\cite{R})}}\label{mincentral}
{Let $X$ be a Hadamard space,     $G$ a group acting    properly and cocompactly by
 isometries on $X$ and $g\in G$ a hyperbolic isometry.  Then $Min(g)$ 
  is  invariant  under  
$C_g(G)$
   and   $C_g(G)$
 acts
   cocompactly on $Min(g)$.}
\end{Th}

     Notice  an orientable irreducible 
closed $3$-manifold    is  a Haken manifold  if it 
contains a  torus or Klein bottle 
  such that the inclusion induces an injective homomorphism between  the fundamental groups.

\begin{Th}\label{3d}
{Let $M$ be a  closed $3$-manifold with nonpositive sectional curvature. Then     every  subgroup 
$H\subset \pi_1(M)$   has the TA-property.}
\end{Th}
\begin{proof}
We may assume $M$ is orientable.  Notice  the theorem follows  from Theorem \ref{haken}  
if    $M$ is a Haken  manifold since
any solvable  subgroup of a group acting   properly and cocompactly  by isometries on a Hadamard space
is virtually free abelian  (see \cite{BH}).
Let $\pi:X\rightarrow M$ be the universal cover of $M$. By a theorem of Eberlein (see \cite{E3}) either $X$  contains a  $2$-flat
  or  $\partial_T X$ is discrete.  If $\partial_T X$ is  discrete, then each   $g\in \pi_1(M)$ is a  
rank one isometry and the theorem follows easily from Corollary \ref{3possibility}.

Suppose $X$  contains a  $2$-flat. Then a theorem of Schroeder (\cite{S2}) says 
  $X$ contains a closed $2$-flat $F$, that is,  $F$ is a $2$-flat and $Stab_{\pi_1(M)}F$ acts cocompactly
  on $F$.  If for any $g\in \pi_1(M)$ either $g(F)=F$ or $g(F)\cap F=\phi$ holds, then
 $\pi(F)$ is an embedded torus or  Klein bottle  and the inclusion induces an injective homomorphism
 on the fundamental groups.
 It follows 
that $M$  is a Haken manifold.

Now suppose  $A$ is a closed   $2$-flat and there is a $g\in  \pi_1(M)$ such that
$g(A)\not=A$ and $g(A)\cap A\not=\phi$. Since  $\dim M=3$  and $A$ and $g(A)$ are totally geodesic, 
  the intersection  $c=A\cap g(A)$ is a complete geodesic. 
 Set $H=Stab_{\pi_1(M)}A$, $K=gHg^{-1}$   and $B=g(A)$.
 Then $K=Stab_{\pi_1(M)}B$  and  $H, K$ and $A, B$ satisfy the assumptions in Theorem \ref{swenson}.
  It follows that $H\cap K$ acts cocompactly on    $c=A\cap B$ and therefore
    must be infinite cyclic.  Let $h$  be  a  generator  of   $H\cap K$ 
 and $P_c$ the parallel set of  $c$. 
 Since $A, B\subset P_c$ and $\dim M=3$, $P_c=Y\times R$  for a surface  $Y$ 
 which is closed   and convex
   in $X$.     $c=\{y_0\}\times R$  for some $y_0\in Y$. 
 $h(P_c)=P_c$  and  $h$  acts  on   $P_c$   as $h=(h_0, t)$, where $h_0$ is an isometry of $Y$
  fixing  $y_0$  and
 $t$ is a translation on $R$. 
Since $h$ leaves $A$  invariant,
    by replacing $h$ with $h^2$ if necessary we may assume  $h_0$ is trivial.  It follows that $P_c=Min(h)$.
 By Theorem \ref{mincentral}   $C_h(\pi_1(M))$ acts cocompactly  on $Min(h)=Y\times R$. 
 If $\p Y=\phi$ then   $X=Y\times R$ and the theorem follows from results on surface groups.

 Suppose $\p Y\not=\phi$.   Let  $Y_0\subset Y$ be the   convex core of $Y$, that is, the smallest 
 closed convex subset of $Y$ with $\p_\i Y_0=\p_\i Y$ (see \cite{L}). Then 
  $Y\subset N_\epsilon(Y_0)$   for some $\epsilon>0$    and $\p Y_0$ consists of
 disjoint complete geodesics.  Note  $C_h(\pi_1(M))$  leaves  $Y_0\times R$  invariant and
 acts on it cocompactly.  Hence there is some $a>0$ such that $d(c_1, c_2)>a$ for any two
 distinct geodesics $c_1$ and $c_2$ in   $\p Y_0$.   Fix a geodesic $c_1\subset \p Y_0$
  and set $F_1=c_1\times R$. Now it suffices to show that
  for any $g\in \pi_1(M)$ either $g(F_1)=F_1$ or $g(F_1)\cap F_1=\phi$.

Suppose there is some $g\in \pi_1(M)$ with  $g(F_1)\not=F_1$   and  $g(F_1)\cap F_1\not=\phi$.
Set $F_2=g(F_1)$.  
Then  $c'=F_1\cap F_2$ is a 
complete geodesic.  If $c'$ is parallel to $c$, 
  then $F_2\subset P_c=Min(h)$, contradicting to 
   the  fact that
 $F_1$, $F_2$ intersect transversally. So $c'$ is not parallel to $c$.
  Choose a geodesic $c''\subset F_2\cap (Y_0\times R)$  
 such that $c''$ is parallel to $c'$ and  $0<d(c', c'')<a$. 
Then $c''\subset \gamma \times R$ for a complete geodesic $\gamma\subset Y_0$.
  Since $c''$   and $c'$  are parallel,   $\gamma$   and  $c_1$  are  parallel and   
 bound a flat strip in $Y_0$.  The fact that $Y_0$ is the convex core of $Y$ implies 
  $\gamma=c_1$ and so $c''\subset F_1$, contradicting to the fact that 
$F_1$, $F_2$ intersect transversally.



\end{proof}

\section{Decomposition of Real  Analytic $4$-Manifolds}\label{deco}

In this section we shall prove   one  of  the main results (Theorem \ref{maint}) of the paper.
The definition of a cycle is given in Definition \ref{cycle}.

\begin{Th}\label{maint}
{Let $M=X/\Gamma$  be a rank 1 closed real analytic 4-manifold of nonpositive sectional curvature.
Suppose there are  no cycles in $X$.  If a subgroup of  $\Ga$   is not virtually free abelian, then it contains a free 
group  of rank two.}
\end{Th}

\begin{remark}\label{exam}
 \emph{There  exist  real analytic   closed    $4$-manifolds (see \cite{AS})
satisfying the   assumptions  of   the theorem.}
\end{remark}

Let $M=X/\Gamma$  be as in Theorem \ref{maint}.
We shall  decompose $X$ into convex domains and construct a 2-complex  $Y$ associated to 
the decomposition.    The 2-complex  $Y$ is a $CAT(-1)$ space and the group $\Gamma$
 acts on $Y$ as a group of isometries.   Thus any subgroup $H$ of $\Gamma$ also
acts on $Y$ as isometries.       A 
    group acting on a $CAT(-1)$ space contains 
a free group of rank two unless it fixes a point in $\ol  Y=Y\cup \partial_\infty Y$ or
stabilizes a geodesic.  Therefore it suffices   to consider these special cases. 

  The decomposition and the associated 2-complex should be considered as 
analogues of  JSJ decomposition and the associated Bass-Serre tree in  $3$-manifold theory.  

\subsection{Decomposition  and the Associated 2-complex}  \label{2complex}

Let $M=X/\Gamma$  be a rank 1 closed real analytic 4-manifold of nonpositive sectional curvature.
In this section we shall decompose $X$,  construct   the associated  2-complex  $Y$ and prove 
$Y$ is a $CAT(-1)$   space.


We  use the results and notation of Section \ref{real}.
Recall 
${\mathcal {W}}$   is  the set of  maximal 
higher rank submanifolds of X that are of the form $Q\times R$ for nonflat 2-dimensional 
Hadamard manifolds $Q$.  
Each $W\in \mathcal {W}$ is called a \emph{wall}. 
We notice 
each    wall   $W$ is totally geodesic and $X- W$ has two components.
Set $X_0=X- \cup_{W\in \mathcal{W}}W$.  We  say  
  two components $C_1$ and $C_2$
  of   $X_0$  are {separated}    by a wall
$W$ if    they    lie in different components of $X- W$, and  say 
$C_1$ and $C_2$   lie on the same side of  $W$ if they lie in the same 
  component of $X- W$. 
Similarly we     define    two 
points $x, y\in X- W$     either   to be   separated by $W$    or  to  lie on the same side of 
  $W$.


Suppose       $W_1$,  $W_2$ are two walls with  
$W_1\cap W_2\not=\phi$. 
  Then    $W_1\cap W_2\cap  W=\phi$ for any  wall 
$W\not=W_1, W_2$.  
   Since each   wall
   is closed, 
  and there are only  a finite number of  walls 
modulo $\Gamma$, there is some $\epsilon>0$ depending
only on $M$ 
       with     $d(W_1\cap  W_2, W)>\epsilon $ for any   wall  
$W\not=W_1,  W_2$.    It   follows that 
   there are 
   4 components $C_1$, $C_2$, $C_3$, $C_4$  of $X_0$    with   
$W_1\cap W_2\subset \ol  C_i$ ($i=1,2,3,4$) and $W_1\cap W_2\cap \ol  C=\phi$ for any
component $C$ of $X_0$, $C\not=C_1, C_2, C_3, C_4$.    By suitably labeling
the 4 components $C_1$, $C_2$, $C_3$, $C_4$, we may assume the following
intersections  are $3$-dimensional:
$\ol  C_1\cap \ol  C_2\subset  W_1$, $\ol  C_2\cap \ol  C_3\subset  W_2$,
$\ol  C_3\cap \ol  C_4\subset  W_1$, $\ol  C_4\cap \ol  C_1\subset  W_2$.

Now we construct the $2$-complex $Y$.        The set $V$ of vertices of $Y$ is in one-to-one
   correspondence    with the set of components of $X_0$.   We denote   by  $C_v$ the component of 
$X_0$ that corresponds to the vertex $v\in V$.     Two   vertices $v_1$ and $v_2$ are joined by 
an edge (denoted   by $v_1v_2$)    if and only if   $\ol  C_{v_1}\cap \ol  C_{v_2}$  is three dimensional.  
 When $W_1$,   $W_2$  are two walls with 
   $W_1\cap W_2\not=\phi$,         there are 
   4 components $C_{v_1}$, $C_{v_2}$,$C_{v_3}$,$C_{v_4}$  of $X_0$ 
such that the following
intersections  are $3$-dimensional:
$\ol  C_{v_1}\cap \ol  C_{v_2}\subset  W_1$, $\ol  C_{v_2}\cap \ol  C_{v_3}\subset  W_2$,
$\ol  C_{v_3}\cap \ol  C_{v_4}\subset  W_1$, $\ol  C_{v_4}\cap \ol  C_{v_1}\subset  W_2$. 
Thus there are the following edges in $Y$:  $v_1v_2$, $v_2v_3$, $v_3v_4$, $v_4v_1$.
 Now whenever there are    two walls  $W_1$,  $W_2$   with
$W_1\cap W_2\not=\phi$,    and $v_1v_2$, $v_2v_3$, $v_3v_4$, $v_4v_1$ the corresponding edges in $Y$ as 
described above,  we attach a square (denoted  by $S(W_1,  W_2)$)  
 along these 4 edges.   The construction of $Y$ is complete.

Next we will put a metric on $Y$.   A \emph{hyperbolic  square} 
  is  a closed convex region  in
the real hyperbolic plane  whose boundary is the union of  4 geodesic segments  
   such that 
 the 4 geodesic segments  are of the same length and the interior angles  at the 
   endpoints of these geodesic segments  are   all  
equal; each of the 4 geodesic segments is called an edge and the endpoints of the 
 4 geodesic segments are called vertices. 
Let $S_0$ be  a  hyperbolic  square so that the interior angle $\alpha$ 
  at the vertices  satisfies:  $\alpha > \frac{2}{5}\pi$.     Let $l$ be the length
of an edge of $S_0$.       Now we declare that all the edges in $Y$ have length
$l$, and all the squares in $Y$ are isometric to  $S_0$.   Thus $Y$ is a piecewise hyperbolic
$2$-complex.  Since there are only a finite number (actually 3 types) of isometry types of cells
in $Y$,  $Y$  with the path metric  is a complete geodesic metric space.   We shall show that
$Y$ is a $CAT(-1)$ space.

    We first look at the intersection $\ol  C_{v_1}\cap \ol  C_{v_2}$ 
   for each edge  $v_1v_2$  of $Y$.

\begin{Le}\label{inside}
{Let $W=Q\times R$   be a wall so that   $W\cap W_1\not=\phi$ for some   wall  $W_1\not=W$.
Suppose  $Z$ is  a component of $W- \cup_{W'\not=W}W'$, where $W'$ varies over all
    walls   that are distinct from $W$.  
Then    the closure $\ol Z$  of $Z$ has the form:  
 $\ol  Z=Q'\times R\subset Q\times R=W$,  
where $Q'\subset Q$ is a closed  
 convex subset of $Q$ and  is the universal cover of a   nonpositively curved 
   compact surface with   closed geodesics on the  boundary.}

\end{Le}

\begin{proof}
Recall (see Theorem \ref{inter}) if $W\cap W'\not=\phi$, then  $W\cap W'=c\times R\subset Q\times R$ for a complete geodesic
  $c$ in $Q$.  And $W\cap W'\cap W''=\phi$   if  $W$, $W'$, $W''$  are  three  distinct   walls. 
Thus $W- \cup_{W'\not=W}W'=Q\times R-\cup_i (c_i\times R)$, where $\{c_i\}$   is a  disjoint collection   of  
 complete geodesics in $Q$.   Since $W$ is closed and there are only a finite number 
of maximal higher rank submanifolds modulo $\Gamma$,   the lemma follows.
\end{proof}

For any wall $W$, let
 $X_W=X- \cup_{W'\not=W} W'$    where $W'$    varies over all
    walls  different from $W$.

\begin{Le}\label{adjacent}
{Let $v_1v_2$ be an edge of $Y$.     Then: \newline
\emph{(i)}  $C_{v_1}$, $C_{v_2}$ are separated by a unique 
wall $W$; \newline
\emph{(ii)}   for any $x_1\in C_{v_1}$ and  $x_2\in C_{v_2}$,  $W$ is the only wall that 
intersects the geodesic segment $x_1x_2$,  and the intersection is transversal;\newline
\emph{(iii)}     if  $E$ is the component of $X_W$
  that contains $C_{v_1}$,   then $C_{v_1}$ and $C_{v_2}$ are the only two components of $X_0$ contained
in $E$.}


\end{Le}

\begin{proof}

(i) Since $v_1$, $v_2$ are distinct, $C_{v_1}$, $C_{v_2}$   are distinct and are thus separated by at least one 
wall.   By the construction of $Y$,  $\ol  C_{v_1}\cap \ol  C_{v_2}$ is three dimensional.  Let $W$ be any wall
that separates $C_{v_1}$  and $C_{v_2}$.  
   Then   $\ol  C_{v_1}\cap \ol  C_{v_2}\subset W^+\cap W^-=W$, where 
$W^+$ and $W^-$ are the closures   of the two components of $X- W$.  Since for any two
    distinct walls 
$W_1$,   $W_2$ the intersection $W_1\cap W_2$ is at most $2$-dimensional,  
 there is exactly one wall $W$ that separates   $C_{v_1}$  and  $C_{v_2}$.   

(ii) Let $W'$    be a wall  with  $W'\cap x_1x_2\not=\phi$.   The intersection must be 
transversal since $W'$ is totally geodesic and $x_1\notin W'$.   Thus $x_1$ and $x_2$ are in 
   different components of $X- W'$.   It follows that $C_{v_1}$ and $C_{v_2}$ are contained 
in  different components of $X- W'$ and $W'$ separates $C_{v_1}$ and $C_{v_2}$. 
Now (i) implies  $W'=W$.  

(iii) From (ii) we see for any  $x_1\in C_{v_1}$ and  $x_2\in C_{v_2}$,  the segment $x_1x_2$ is contained 
in $E$.  Thus $E$ contains both $C_{v_1}$ and $C_{v_2}$.   Let $C\not=C_{v_1}$ be a component of    $X_0$
contained in $E$.   $C\not=C_{v_1}$  implies $C$ and $C_{v_1}$ are separated by a wall, which must be 
$W$ by the definition of $E$.  Now there is no wall separating $C_{v_2}$ and $C$ and we 
  have $C=C_{v_2}$.

\end{proof}

By the   Cartan-Hadamard theorem, to prove $Y$  is $CAT(-1)$ it suffices to show 
$Y$ is simply connected   and 
has curvature 
$\le -1$.    Since $Y$ is piecewise hyperbolic,  $Y$ has curvature 
$\le -1$ if all the vertex links are $CAT(1)$.   The vertex links are metric graphs  and a 
metric graph is $CAT(1)$ if and only if each   injective  edge  loop has length at least $2\pi$.   
   By the construction of $Y$,  the  edges in the vertex  links all have the same length
 $\alpha>\frac{2}{5}\pi$.       Therefore the vertex links are $CAT(1)$ if each  injective  edge  loop 
(in the links)    has at least 5  edges.    

\begin{Le}\label{link}
{Let $v\in Y$ be a vertex of $Y$ and $Link(v,Y)$ the link of $v$ in $Y$. 
Then each  injective  edge  loop in   $Link(v,Y)$
     has at least 5  edges.    
  In particular, $Y$ has curvature $\le -1$.}
\end{Le}

\begin{proof}
Let $e$ be an edge in $Link(v,Y)$.  Then $e$ corresponds to a square $S$ of $Y$ that 
has $v$ as one of its 4 vertices.  Let $v$, $v_1$, $v'$, $v_2$ be the 4 vertices 
     of   $S$ in cyclic order.   By Lemma \ref{adjacent},  there is a unique
wall $W_1$ that separates $C_v$ and $C_{v_1}$ and a unique wall $W_2$
that separates $C_v$ and $C_{v_2}$.   By the construction of $Y$, $W_1\cap W_2\not=\phi$
   and   the square  $S$  is   determined by  $W_1$   and  $W_2$.    
   The two endpoints of $e$ uniquely determine $v_1$ and $v_2$, 
$v_1$ and $v_2$   then uniquely determine 
$W_1$ and $W_2$, and $W_1$ and $W_2$ in turn uniquely determine the square  $S$.    
   It follows that the edge  $e$  is uniquely
determined by   its  two endpoints. 
Therefore there is no  injective edge loop   with length 2 in $Link(v,Y)$.

Now   let $L$   be    an injective  edge loop   in $Link(v,Y)$ consisting of $n$ edges.
    Then there  are $n$ edges $vv_i$ ($i=1,2, \cdots, n$) in $Y$ so that 
 $vv_i$, $vv_{i+1}$ (indices are taken modulo $n$) are  two edges of a square.
It follows from  the preceding  paragraph that there are 
walls $W_1, \cdots, W_n$      such that 
$W_i$ is  the only wall separating $C_v$ and $C_{v_i}$   and 
   $W_i\cap W_{i+1}\not=\phi$ for all $i$ 
where indices are taken modulo $n$.    Suppose there are indices $i$ and $j$   with  
   $i\not=j$ and $W_i=W_j$.   
     Let   $E$    be  the component of $X_{W_i}$
   containing  $C_v$.   Lemma \ref{adjacent}  applied to    the edge $vv_i$ implies 
  $C_v, C_{v_i}\subset E$. Similarly       $C_v, C_{v_j}\subset E$  as  $W_i=W_j$. 
  Since $L$ is an injective loop in $Link(v,Y)$,  $i\not=j$   implies  $v_i\not=v_j$.
  Therefore  $C_{v_i}\not=C_{v_j}$ and there are three distinct components 
  $C_v, C_{v_i}, C_{v_j}$   of $X_0$ contained in $E$,  contradicting to Lemma \ref{adjacent}.
   It  follows that $W_1, \cdots, W_n$   are all distinct and   form    
      an   $n$-cycle in $X$. Now 
Proposition \ref{no4cycle} completes the proof.

\end{proof}



\begin{Le}\label{1connected}
{Any   edge  loop in $Y$   is homotopically  trivial.
  In particular,   $Y$ is simply connected.}

\end{Le}

\begin{proof}We induct on the length of   an     edge loop. 
   By the construction of $Y$, two vertices are connected by
at most one edge and there is no injective
edge loop with length 2.

Let $l$ be an  injective edge loop in $Y$ and $v_1,  \cdots, v_n$ the vertices 
on $l$ in cyclic order.    Let $x_i\in  C_{v_i}$   be 
an arbitrary  point.
   Notice the components $C_{v_i}$ ($i=1,  \cdots, n$)  are all distinct since 
$l$  is  an  injective edge loop. 
  Denote    by   $W_i$   the unique wall that 
  separates $C_{v_i}$ and $C_{v_{i+1}}$ (here indices are taken
modulo $n$).
  By Lemma \ref{adjacent}  $W_i$ is the only   wall intersecting  $x_ix_{i+1}$.  
  Let  $L=\cup_{i=1}^n x_ix_{i+1}$.    Clearly 
  $L$ is a loop in $X$.

  Suppose   $n=3$.   Since  $W_1$ separates $x_1$ and $x_2$, the path $x_2x_3* x_3x_1$ must
intersect $W_1$.   From the preceding paragraph  the path $x_2x_3* x_3x_1$  only 
   intersects $W_2$ and $W_3$, 
we have $W_2=W_1$ or $W_3=W_1$.   We may assume $W_2=W_1$,  the other case being handled similarly.   Let
$E_2$ be  the component of $X_{W_2}$
  that contains $C_{v_2}$.   Then Lemma \ref{adjacent} implies  the three distinct components 
   $C_{v_1}$, $C_{v_2}$, $C_{v_3}$  are all
  contained in   $E_2$ since $W_2=W_1$,   which contradicts to the same lemma.
 It follows that   there is no    injective  edge   loop with length 3.

Now  suppose  $n\ge 4$.   
Since $W_1$ separates $x_1$ and $x_2$, the path $L-x_1x_2$ must
intersect $W_1$.    By construction the path $L- x_1x_2$
  only intersects   $W_2, \cdots, W_n$.  So $W_i=W_1$ for at least one  $i$,  $2\le i\le n$.
  Let $m$ be the largest such $i$.    
Then $x_1$ and $x_{m+1}$  lie   on   one  side  of 
$W_1$, and $x_2$ and $x_m$ lie     on   the other   side of $W_1$.
The argument in  the preceding paragraph  shows $x_1\not= x_{m+1}$ and $x_2\not=x_m$ since $l$ is an injective loop.
Let $z=x_1x_2\cap W_1$      and   $z'=x_mx_{m+1}\cap W_1$.   Since   $x_2\not=x_m$, we have 
$C_{v_m}\not=C_{v_2}$,  and   by   Lemma \ref{adjacent} (iii)   there is at least one  wall 
$W'$, $W'\not=W_1$ with $W'\cap zz'\not=\phi$.
  Let $\{W'_1,  \cdots, W'_k\}$ be the set  of walls different from $W_1$ that intersect
$zz'$.    Then  $W_1\cap W_i'\not=\phi$ ($i=1, \cdots, k$) and $S(W_1,   W_i')$ is 
  a square in $Y$. 
Set $z_i=zz'\cap W'_i$.  We 
  label the walls $W'_1$,  $\cdots$, $W'_k$  so that $d(z, z_i)<d(z, z_j)$    whenever $i<j$.  
   Then the interior of $zz_1$  does not intersect any wall different from $W_1$. 
  It  follows that  $v_1$, $v_2$     are  vertices  of  $S(W_1,   W_1')$.
Let $v'_2$,   $v''_2$ be the other two vertices of $S(W_1,   W_1')$   
so that  $C_{v'_2}$,  $C_{v_2}$
lie on the same  side of $W_1$.  Similarly we define $v'_i$, $v''_i$  for all
 $2\le i\le k+1$   so that  $v'_{k+1}=v_m$, $v''_{k+1}=v_{m+1}$  and 
 all $C_{v'_i}$ ($2\le i\le k+1$)
  lie on the same side of $W_1$.

  Let
$l_1=v_2v'_2* \cdots * v'_kv_m$,  $l_2=v_1v''_2* \cdots * v''_kv_{m+1}$,
$l_3=v_2v_3* \cdots *  v_{m-1}v_m$, $l_4=v_{m+1}v_{m+2}* \cdots * v_nv_1$
 be oriented paths in $Y$ and $l^{-1}_i$ ($i=1,2,3,4$)   the same paths with the reverse orientation.
Set $l'=v_1v_2* l_3* l^{-1}_1* v_2v_1$, 
$l''=v_1v_2* l_1* v_mv_{m+1}* l^{-1}_2$,
$l'''=l_2*   l_4$.   Then   $l'$, $l''$ and $l'''$ are oriented loops so that 
$l=v_1v_2* l_3* v_mv_{m+1}* l_4$ is   homotopic to $l'* l''* l'''$.
 Notice  $\cup_{i=1}^k S(W_1, W'_i)$  is homeomorphic to a square and $l''$ is its boundary.
    Therefore $l''$ is homotopically trivial.   
  Now we notice $length(l_3)\ge k$ since the   path $x_2x_3* \cdots x_{m-1}x_m$
  must cross all the walls $W'_1, \cdots, W'_k$.  Similarly $length(l_4)\ge k$.
Since $length(l_1)=length(l_2)=k$,  the lengths of the loops  $l_3* l^{-1}_1$
  and $l_2* l_4$ are strictly less than the length of $l$. The induction hypothesis 
  implies   $l_3* l^{-1}_1$
  and $l_2* l_4$ are homotopically trivial. Therefore $l'$ and $l'''$ are 
homotopically trivial. It follows $l$ is also homotopically trivial.

\end{proof}

Lemmas \ref{link} and \ref{1connected} together imply the following proposition.

\begin{Prop}\label{cat-1}
\emph{ The $2$-complex $Y$ is a $CAT(-1)$ space.}
\end{Prop}

Recall a group   is   said to have  the TA-property if it is virtually 
free abelian or contains a free group of rank two.   
We note Theorem \ref{maint} holds if and only if any subgroup $H$ of $\Gamma$ has the TA-property.
   Let $H\subset \Gamma$  be a  subgroup of $\Gamma$. 
By Lemma \ref{finiteid} 
 we may assume    each $h\in H$ preserves the orientation of $X$,
after  passing to an  index two subgroup if necessary.
    Corollary \ref{1possibility} implies we may assume $H$  does not contain 
  any rank one  isometry. By Theorem \ref{kim}, we may further assume there is a
Tits component $\mathcal C\subset \partial_T X$ so that $h(+\infty), h(-\infty)\in {\mathcal C}$ for all
$h\in H$.   If ${\mathcal C}=\partial_T F$ for a maximal higher rank submanifold    $F$  which 
is a  $2$-flat or $3$-flat,   then each element $h\in H$ leaves  $F$  invariant. 
  It follows that $H$ acts on the Euclidean space $F$ properly and isometrically, and thus
must be virtually free abelian by Bieberbach's theorem.




From now on we assume $H\subset \Gamma$ satisfies the following properties:\newline
(a)   each  $h\in H$ preserves the orientation of $X$;\newline
(b) $H$  does not contain 
  any rank one  isometry;\newline
(c)  there is a  unique 
Tits component $\mathcal C\subset \partial_T X$    such that   $h(+\infty), h(-\infty)\in {\mathcal C}$ for all
$h\in H$;\newline
(d)  the Tits component $\mathcal C$ has the following form: 
 $\mathcal C=\bigcup_{W\in {\mathcal W}^*}\partial_{\infty} W$,  where $ {\mathcal W^*}\subset {\mathcal {W}}$ 
  is a subset such that  $\bigcup_{W\in {\mathcal W}^*}W$ is a connected component of 
$\bigcup_{W\in {\mathcal W}}W\subset X$.

We shall  also   frequently pass to finite index subgroups of $H$.

\subsection{Action on the $2$-complex}\label{actony}

Since $\mathcal{W}$ is invariant under the action  of  the group $\Gamma$,  it is
clear from the construction of $Y$ that $\Gamma$ acts on $Y$ as a group of 
  cellular  isometries. 
Any  cellular isometry of $Y$   either     is   hyperbolic or have a fixed point in $Y$ as 
$Y$ only has a finite number of isometry types of cells (see \cite{Br}).
  The fact  that  $Y$ is a $CAT(-1)$ space implies any hyperbolic isometry
of $Y$ is of rank one.  Thus  any $\gamma\in \Gamma$ either acts on $Y$ as 
a rank one isometry or has a fixed point in  $Y$.     Notice 
an element  $\gamma\in \Gamma$    may
act  as a rank one isometry on $Y$
   even if  it is not rank one  (with respect to its action on
  $X$) in $\Gamma$.

 For any $\gamma\in \Gamma$,
  let $\gamma_Y:Y\rightarrow Y$ denote the  induced isometry of $\gamma$ on $Y$, and 
$Fix(\gamma_Y)\subset Y$ denote the fixed point set of $\gamma_Y$. 
Let $H\subset \Gamma$  be a subgroup    satisfying  the properties stated  
     at the end  of   Section \ref{2complex}.
As a  subgroup of $\Gamma$, $H$ also acts on $Y$ as isometries. 
We shall   analyze
the action of $H$ and  its individual elements on $Y$.


\begin{Le}\label{square}
{Let    
$S=S(W_1, W_2)$    be  a square in $Y$  
and  
   $\gamma\in \Gamma$.   
  If  $\gamma_Y(S)=S$, then $\gamma(W_1 \cap W_2)=W_1\cap W_2$.}
\end{Le}
\begin{proof}
Since   $W_1,  W_2$  are uniquely determined by  $S$,   $\gamma_Y(S)=S$ implies 
 $\{\gamma(W_1), \gamma(W_2)\} =\{W_1, W_2\}$. The lemma follows.

\end{proof}

Let $e=v_1v_2$ be an edge in $Y$, and $C_1$, $C_2$ the two components of $X_0$
 corresponding to $v_1$, $v_2$ respectively.  Lemma \ref{adjacent}  implies there is a unique
wall separating $C_1$ and $C_2$.   The following  lemma is clear.

\begin{Le}\label{edge}
{Let $e=v_1v_2$ be an edge in $Y$,  $C_1$, $C_2$ the two components of $X_0$
 corresponding to $v_1$, $v_2$ respectively, and    $W$  the  unique
wall separating $C_1$ and $C_2$.  If $\gamma_Y(e)=e$ for some $\gamma\in \Gamma$, then 
$\gamma(W)=W$.}
\end{Le}

   Suppose    $\gamma\in \Gamma $  is not a rank one isometry. 
  Let $P(\gamma)$ denote the parallel set of  an axis of $\gamma$.  
Since $X$ is real analytic, $P(\gamma)$ is a higher rank submanifold
in $X$, and therefore is contained in a maximal higher rank submanifold. 
We further suppose $P(\gamma)\subset W$ for a wall 
 $W=Q\times R$, where $Q$ is 
 a $2$-dimensional nonflat Hadamard manifold.  
   Notice $\gamma$  belongs to  exactly
  one of the following  4    classes:\newline
Type A: $P(\gamma)=c\times R\subset Q\times R=W$ and $P(\gamma)\cap W'=\phi$ for any 
   wall $W'\not=W$, 
 where $c$ is a complete geodesic in $Q$;\newline
Type B: $P(\gamma)=W\cap W'$ for some  wall $W'\not=W$;\newline
Type C: $P(\gamma)=W$. In this case,   $\{q_0\}\times R$
 is an axis of  $\gamma$   for  some  $q_0\in Q$;\newline
Type D: $P(\gamma)=c\times R\subset Q\times R=W$ and    the axes  of $\gamma$    intersect transversally with
 some  wall $W'$,   where $c$ is a complete geodesic in $Q$.   


Before we study the action of  individual isometries on $Y$, we introduce a subcomplex 
of $Y$ associated to each wall that intersects other walls.   
Recall     for each   edge   $v_1v_2$  in $Y$ there is a  unique wall $\tilde W$ separating 
 $C_{v_1}$ and $C_{v_2}$.  
  Given any square $S(W_1, W_2)$ in $Y$, 
 two opposite edges of $S(W_1, W_2)$  determine $W_1$, while the other two opposite edges
determine $W_2$.   Let $\omega(W_1,W_2)\subset S(W_1,W_2)$
 be the geodesic segment 
  connecting the midpoints of the two opposite edges of $S(W_1, W_2)$ 
that determine $W_1$.  
Let  $W$   be   a wall intersecting   other walls.
   Set     $C_W=\bigcup S(W,W')$   and $T_W=\bigcup\omega(W,W')\subset C_W$ 
 where $W'$    varies over all
walls $W'\not=W$
with  $W'\cap W\not=\phi$.    $C_W$ is   a subcomplex of $Y$.


\begin{Le}\label{sw}
{Let   $\gamma\in \Gamma$  be  a Type C 
   isometry.  If  the wall    $W:=P(\gamma)$
 intersects other walls,      then there is an integer $n$ such that 
$Min(\gamma^n)=W$   and 
   $Fix(\gamma^n_Y)=C_W$.    In particular
$C_W$   is a  closed  convex subset of $Y$, and  $\phi\not=Fix(\gamma_Y)\subset C_W$.}
\end{Le}
\begin{proof}
Let $W=Q\times R$. 
By the definition of a  Type C isometry, there is some 
$q_0\in Q$
  such that $\{q_0\}\times R$
 is an axis of  $\gamma$.  
 It follows that $\gamma_{|W}$ has the  following form 
$\gamma_{|W}=(e,t):Q\times R\rightarrow Q\times R$  where $e$ is an elliptic isometry of $Q$ fixing
the point $q_0$ and $t$ is a translation of $R$.  
Since the  wall $W$ is closed  and $Q$ is a
nonflat $2$-dimensional
 Hadamard manifold,   $e$ has finite order (see \cite{E1}).  Therefore there is an even integer
  $n$ such that $e^n=id$. Notice $Min(\gamma^n)=W$  and $\gamma^n$ preserves the 
orientation of $X$.

    Now assume   $\gamma\in \Gamma$  is    a Type C 
   isometry  preserving  the 
orientation of $X$  
with $Min(\gamma)=W$.  We shall show $Fix(\gamma_Y)=C_W$.
We first argue $C_W\subset Fix(\gamma_Y)$.  Let $W'\not=W$ with $W'\cap W\not=\phi$. 
Since $Min(\gamma)=W=Q\times R$, each geodesic  $\{q\}\times R$ ($q\in Q$) is  an axis of $\gamma$. 
 It follows that the  $2$-flat $W\cap W'=c\times R\subset Q\times R=W$ is invariant under $\gamma$. 
  Since  $W'$ is perpendicular to $W$, we have $\gamma(W')=W'$.
          Therefore   $\gamma_Y(S(W,W'))=S(W,W')$.  
  Since $\gamma$  
 preserves the orientation of $X$ and translates all the geodesics $\{q\}\times R$ ($q\in Q$),
   each of the 4 components of $X_0$ determined by $W$ and $W'$ is invariant under $\gamma$.
So the square $S(W,W')$ is  pointwise   fixed   by $\gamma_Y$.  

Next  we show  $ Fix(\gamma_Y)\subset C_W$.  Let $y\in Fix(\gamma_Y)$.  If $y$ lies in the interior of 
   a square $S(W_1, W_2)$, then $\gamma_Y(S(W_1,  W_2))=S(W_1, W_2)$.  
 Lemma \ref{square} implies $\gamma(W_1\cap W_2)=W_1\cap W_2$ and  thus an axis of $\gamma$ lies in 
$W_1\cap W_2$.  It follows that $W\cap W_1\cap W_2\not=\phi$.  By Theorem \ref{subm} (iv)
   $W=W_1$ or $W=W_2$.  In either case $ S(W_1, W_2)\subset C_W$.

 Suppose  $y$ lies in the interior of an edge $v_1v_2$. Then $\gamma_Y(v_1v_2)=v_1v_2$. 
   Let  $W'$   be the unique wall 
  separating $C_{v_1}$ and $C_{v_2}$.   Lemma \ref{edge} implies $\gamma(W')=W'$ and thus $W'$ contains 
an axis of $\gamma$.  It  follows  $W\cap W'\not=\phi$.  
 If $W'=W$, then clearly $v_1v_2\subset C_W$. Suppose $W'\not=W$.  Pick
$x_1\in C_{v_1}$, $x_2\in C_{v_2}$ and let $z'=x_1x_2\cap W'$. Also pick $z\in W'\cap W$.  
Then $zz'\subset W'$.  Let $W=W_1$, $\cdots$, $W_k$ be the sequence of walls  in consecutive order 
    that intersect 
$zz'$ transversally.  These walls   determine a sequence of squares   
$S(W',W_1)$, $\cdots$,  $S(W',W_k)$  in $Y$, where two   consecutive  squares
share exactly an edge.  Notice   $v_1v_2$ is an   edge of $S(W',W_k)$.
 For each of these squares, two opposite edges
correspond to the wall $W'$.  Let $l_i$ ($i=1, 2,\cdots, k$) be the geodesic segment connecting the midpoints 
of  these two opposite edges,  and $l=\cup_{i=1}^k l_i$. 
$l$ is a geodesic  by the geometry of $Y$. One endpoint   of
  $l$ lies in $S(W',W)\subset C_W$ and the other one is the midpoint of the edge $v_1v_2$, so they  are both
fixed by $\gamma_Y$. It follows that $l$ is pointwise fixed by $\gamma_Y$, and by the preceding 
    paragraph    $S(W',W_i)\subset C_W$ for each $i$. In particular  $v_1v_2\subset C_W$.

  Now suppose $y$ is a vertex. Then for any $y'\in C_W$, the  geodesic segment $yy'$ is pointwise  fixed 
by $\gamma_Y$.  The initial segment of $yy'$ in contained in some square or edge. It 
  follows from the above  that 
this square or edge lies in $C_W$ and so does $y$.    

\end{proof}

It follows from Lemma \ref{sw} that  $C_W$  is a $CAT(-1)$ space.  $C_W$ clearly admits 
a  \lq\lq reflection" about $T_W$ with $T_W$ as the fixed point set.  Thus $T_W$ is also a closed and convex subset of $Y$.
$T_W$ is the  \lq\lq core" of $C_W$, and $C_W$ is homeomorphic to $T_W\times [0,1]$.  For 
each wall  $W$  that intersects  other walls, 
  let $T_W(\infty)\subset \partial_\infty Y$ be the  geometric  boundary of  $T_W$ 
     naturally identified  with a subset
   of $ \partial_\infty Y$.  

\begin{Le}\label{swinter}
{Let $W_1\not=W_2$ be two  walls  that intersect other walls. Then  
$T_{W_1}(\infty)\cap T_{W_2}(\infty)=\phi$.}
\end{Le}

\begin{proof}
Suppose $T_{W_1}(\infty)\cap T_{W_2}(\infty)\not=\phi$   and pick
  $\xi\in T_{W_1}(\infty)\cap T_{W_2}(\infty)$.     Let
$\alpha_1: [0,\infty) \rightarrow T_{W_1}$ and $\alpha_2: [0,\infty) \rightarrow T_{W_2}$
  be two rays ending at $\xi$.   Since $Y$ is a $CAT(-1)$ space, 
we have   $d(\alpha_1(t), \alpha_2(t))\rightarrow 0$ as $t\rightarrow \infty$.    
Since  $\alpha_1$  and $\alpha_2$   cross  centers of squares    in $Y$,   there are  $a>0$, $b\in R$
    so that
   $\alpha_1(t)=\alpha_2(t+b)$ for  all $t\ge a$.    $\alpha_1$    passes through midpoints 
of edges in $Y$ and   such a  midpoint  uniquely
determines the wall $W_1$.   The same is true for  $\alpha_2$   and $W_2$.   
 Hence $W_1=W_2$,  a contradiction.

\end{proof}

  Let  $\gamma$   be   a Type A isometry and $W$  the  wall with $P(\gamma)\subset W$. 
   Then  $P(\gamma)\subset X_W$. Let $E$ be the component of 
$X_W$ containing  $P(\gamma)$, and   $C_1$ and $C_2$    the two components of  
$X_0$ contained in $E$.   Denote    by   $v_1$ and $v_2$ the two vertices of $Y$ 
corresponding to  $C_1$ and $C_2$  respectively.     Then  $v_1v_2$ is an edge in $Y$.

\begin{Le}\label{typea}
{ Let $\gamma$ be a Type A isometry, and $W$, $v_1v_2$ be as above.
Then    $\phi\not=Fix(\gamma_Y)\subset v_1v_2$.}
\end{Le}
\begin{proof}
We use the notation from  the paragraph preceding the lemma.
Since $E$ is  the component of 
$X_W$ containing  $P(\gamma)$   and $P(\gamma)$ is invariant under 
$\gamma$,  $\gamma(E)=E$.     It follows  that  $\{\gamma(C_1), \gamma(C_2)\}=\{C_1, C_2\}$
  and $\gamma_Y$ sends $v_1v_2$ to  itself.  In particular, the midpoint of $v_1v_2$ is
  fixed by $\gamma_Y$.   Lemma \ref{square} implies   no square in $Y$ is invariant
under   $\gamma_Y$.    Then  Lemma  \ref{edge}   and an 
    argument similar to   the proof  
of Lemma \ref{sw} shows $Fix(\gamma_Y)\subset v_1v_2$. 

\end{proof}

\begin{Le}\label{typeb}
{Let $\gamma$ be a Type B isometry and $W_1$, $W_2$ the two walls   such that 
$P(\gamma)=W_1\cap W_2$.   Then $\phi\not=Fix(\gamma_Y)\subset S(W_1, W_2)$.}
\end{Le}
\begin{proof}
Since $P(\gamma)=W_1\cap W_2$,  we have  $\{\gamma(W_1), \gamma(W_2)\}=\{W_1, W_2\}$
    so  the square $S(W_1, W_2)$ is invariant under $\gamma_Y$. In particular, 
the center of the square  $S(W_1, W_2)$  is fixed by $\gamma_Y$.  
  Now an  argument  similar to the proof of    Lemma \ref{sw}
   shows   $ Fix(\gamma_Y)$   is contained in $S(W_1, W_2)$.

\end{proof}

\begin{Le}\label{typed}
{Let $\gamma$ be a Type D isometry  with $P(\gamma)$ contained in  a wall $W$.
 Then $\gamma_Y$ is a rank one isometry of $Y$,
and the axis of $\gamma_Y$ is contained in $T_W$.}

\end{Le}

\begin{proof}
 We exhibit a geodesic in $T_W$ that is translated by $\gamma_Y$.  
Let  $c\subset P(\gamma)\subset  W$ be  an axis of $\gamma$.
Then $c$ intersects some wall $W'$
transversally   since $\gamma$ is a Type D isometry.
Pick a point $p\in c$ such that $W$ is the only wall containing $p$. 
 Let $W_1$, $W_2$, $\cdots$, $W_k$ be the sequence of walls in consecutive order 
 intersected transversally 
 by the geodesic segment $p\gamma(p)$.  Then the biinfinite  sequence  of walls
$$\cdots,\gamma^{-1}(W_1), \cdots, \gamma^{-1}(W_k), W_1, \cdots, W_k, \gamma(W_1), \cdots, \gamma(W_k), \cdots$$
 is the sequence of walls   in consecutive order  
intersected   transversally by the   axis  $c$. This sequence of walls   together with $W$ determine
    a sequence of squares in $Y$. The union of this sequence of squares is a \lq\lq strip"  $R$
  contained in $C_W$, and $\gamma_Y(R)=R$.   It follows that 
the geodesic  $R\cap T_W$, the center line of $R$,   is translated by  $\gamma_Y$.  

\end{proof}

Let $H\subset \Gamma$  be a subgroup  satisfying the properties stated   at the end  of   Section \ref{2complex}.
As a  subgroup of $\Gamma$, $H$ also acts on $Y$ as isometries.  
  $Y$ is a Gromov hyperbolic space   since  it is  $CAT(-1)$.  By a theorem of Gromov (see \cite{G}),   one of the following 
occurs:\newline
(1) $H$ has a bounded orbit in $Y$;\newline
(2) $H$ has a fixed point in $\partial_\infty Y$;\newline
(3) $h_Y$ is a rank   one isometry  for  some $h\in H$.

When (1) occurs,  Cartan's fixed point theorem implies $H$ has a fixed point in $Y$. 
After passing to a  finite index subgroup, we may assume that $H$ fixes  a vertex of $Y$. 
In Section \ref{fixedpoint} we show  $H$ has the TA-property  if it fixes a vertex of $Y$
  and there is no cycle in $X$.

\begin{Le}\label{v}
{ Suppose $H$ has a fixed point in $\partial_\infty Y$ and  $h_Y$ is   not  a rank one isometry
 for    any   $h\in H$. Then $H$   has the  TA property.}

\end{Le}
\begin{proof}
By assumption and the analysis of the 4 types of isometries,   each  $h\in H$  is of Type C.
  By Lemma \ref{sw}  if $h$ is of Type C  and $W$ is the wall with $P(h)=W$, 
then the fixed point  set of $h_Y$  in $\partial_\infty Y$     
  is  contained in   $T_W(\infty)$.   Now Lemma \ref{swinter} and the fact that $H$ has a fixed point in 
$\partial_\infty Y$  imply that there  is a wall $W$  with  $P(h)=W$ for all $h\in H$.  
  Therefore    $h(W)=W$  for all   $h\in H$ 
   and  $H\subset Stab_{\Gamma}W$.  Since $Stab_{\Gamma}W$
  is the fundamental group of the   closed $3$-manifold  $W/Stab_{\Gamma}W$
 with nonpositive  sectional curvature,  the lemma follows from Theorem \ref{3d}.

\end{proof}

Now we remain to consider  the case when $g_Y$ is rank one for some $g\in H$.  
Let $c\subset Y$ be the axis of $g_Y$, and $c(+\infty)$,  $c(-\infty)$  
  the fixed points  of $g_Y$ in $\partial_\infty Y$.
By    Corollary \ref{3possibility} one of the 
 following holds:\newline
(1) $c(+\infty)$   or   $c(-\infty)$  is fixed by all elements of $H$; \newline
(2)  $h_Y(c)=c$ for all $h\in H$; \newline
(3)   $H$ contains  a free group of rank two.

We only need to consider the first two cases.    Case (2) can be reduced to case (1)
   as  in case (2) an index two subgroup of $H$ fixes   both  
$c(+\infty)$ and $c(-\infty)$.

\begin{Le}\label{atinf}
{Suppose $g_Y$ is rank one for some $g\in H$,  and $c(+\infty)$,  $c(-\infty)$  
   are   the fixed points  of $g_Y$ in $\partial_\infty Y$.   If $c(+\infty)$  
 or   $c(-\infty)$  is fixed by all elements of $H$, then   $H$   has the TA-property.}

\end{Le}
\begin{proof}
The assumption implies $g$ is of Type D.   Let $W$ be the wall with 
  $Min(g)\subset W$.  Then the axis $c$ of $g_Y$ is contained in $T_W$.   
 Since each $h\in H$ fixes a point in $\partial_\infty Y$,  
by the analysis of the 4 types of isometries we see each $h\in H$ is either of Type D or 
  of Type C.    By Lemmas \ref{sw}, \ref{swinter}, \ref{typed} and the fact that $c(+\infty)$  
 or   $c(-\infty)$  is fixed by all elements of $H$,     we conclude
 $h(W)=W$ for all $h\in H$.   Now the  argument in the proof of Lemma \ref{v}
 shows    $H$   has the TA-property.

\end{proof}

\subsection{Subgroups With a Fixed Point in the $2$-complex}\label{fixedpoint}

We remain to consider the case when $H$ fixes a vertex of $Y$. To be  more precise,
  throughout this section, $H\subset \Gamma$  is a subgroup 
satisfying the following properties:\newline
(a)  there is a  
Tits component $\mathcal C\subset\partial_T X$ so that $h(+\infty), h(-\infty)\in {\mathcal C}$ for all
$h\in H$;\newline
(b)  $\mathcal C=\bigcup_{W\in {\mathcal W}^*}\partial_{\infty} W$,  where $ {\mathcal W^*}\subset {\mathcal {W}}$ 
  is a subset such that  $\bigcup_{W\in {\mathcal W}^*}W$ is a connected component of 
$\bigcup_{W\in {\mathcal W}}W\subset X$;\newline
(c)  there is a vertex $v$ in $Y$ so that $h_Y(v)=v$ for all $h\in H$.

 Notice condition (c) implies $H$ contains no Type D isometries since by Lemma \ref{typed}
   $\gamma_Y$  has no fixed point in $Y$ if $\gamma$ is of Type D.

 Set $W_{\mathcal C}=\bigcup_{W\in {\mathcal W^*}}W$.
        Then   both $C_v$ and  $W_{\mathcal C}$ are 
invariant under $H$. It follows that $\ol  C_v\cap W_{\mathcal C}$ 
  is invariant under $H$.  Set $\Lambda=\ol  C_v\cap W_{\mathcal C}$. Then $H\subset Stab_{\Gamma}\Lambda$. 
    We  shall prove  $\Lambda$   is connected    and       $\Lambda/{Stab_{\Gamma}\Lambda}$
   is a compact  $3$-manifold.    Under the assumption that there is no
  cycle in $X$,   $\Lambda$  is simply connected and consequently
  $H$ is a subgroup of  the fundamental group of a compact $3$-manifold. 
    Theorem  \ref{haken}    then  implies    $H$ 
     has the TA-property.

For any wall $W$, the closures of the two components of $X-W$ are called closed half spaces.
  $\ol C_v$ is clearly the intersection of a family of closed half spaces.
Now it   is  not hard to see that if  $W$   is  a wall with  $W\cap \ol  C_v\not=\phi$, then
    $W\cap \ol  C_v=\ol  Z$  for   a component    $Z$ of $W- \bigcup_{W'\not=W}W'$,
   where $W'$  varies  over all walls    distinct from $W$.  By Lemma \ref{inside}  such a $\ol  Z$ 
  has the form $\ol  Z=Q'\times R\subset Q\times R=W$,  where $Q'\subset Q$ is a closed 
 convex subset of $Q$ and  is the universal cover of a  nonpositively curved 
   compact surface with   closed geodesics on the  boundary. 
    Similarly    we  see  if  $W_1\not=W_2$   and  
    $W_1\cap W_2\cap \ol  C_v\not=\phi$,   then 
 $W_1\cap W_2\cap \ol  C_v$   is a $2$-flat.

Let $W, W'\in   {\mathcal W^*}$.
A   \emph{chain}  from $W$ to $W'$   is a sequence 
   $W=W_0$, $W_1$, $\cdots$, $W_n=W'$    of  walls   in $  {\mathcal W^*}$   
with $W_i\cap W_{i+1}\not=\phi$ ($i=0, \cdots, n-1$).   Since $W_{\mathcal C}$  is  connected,
   there is   at   least  one chain from  $W$ to $W'$.

\begin{Le}\label{connected}
{  $\Lambda$  is path connected.}


\end{Le}
\begin{proof} 
Let $W, W'\in   {\mathcal W^*}$ with  $W\cap \ol  C_v, W'\cap \ol  C_v\not=\phi$,   and 
   $W=W_0$, $W_1$, $\cdots$, $W_n=W'$   a chain from $W$ to $W'$  with the smallest possible $n$.  
  Choose  $x_0\in W\cap \ol  C_v$ and $x_{n+1}\in W'\cap \ol  C_v$  so that
   $x_0$ does not lie in any wall other than $W$ and $x_{n+1}$ does not
 lie in any wall other than $W'$.
  Pick $x_i\in W_{i-1}\cap W_i$  ($i=1,\cdots, n$)    and    let
  $\sigma$   be  the path  defined by 
$\sigma=x_0x_1* \cdots* x_nx_{n+1}$.
  Note the lemma follows if 
$\sigma\subset \ol  C_v$.

We   show  $\sigma\subset \ol  C_v$ by  inducting  on the length of  a  chain. 
We first show $x_0x_1\subset \ol  C_v$.  It suffices to show $x_1\in \ol  C_v$ since $\ol  C_v$ is convex.
Assume $x_1\notin \ol  C_v$. 
Then there is a  wall $\tilde W\not=W, W_1$  so that $x_0$ and $x_1$  lie in 
   different components of  $X-\tilde W$.
  Since $x_0x_1\subset W$,  $W\cap \tilde W\not=\phi$. 
  $\tilde W\not=W'$  holds, otherwise by the choice of $n$ we have $\tilde W=W_1$,
     a contradiction.   It follows that  $x_0$ and $x_{n+1}$ lie on  the same 
  side   of $ \tilde W$.  Thus the part of $\sigma$ from $x_1$ to $x_{n+1}$
  must cross $\tilde W$.  Suppose $\tilde W\cap x_ix_{i+1}\not=\phi$ for some $i$, $1\le i\le n$.  
Since  $x_ix_{i+1}\subset W_i$ we have $\tilde W\cap W_i\not=\phi$. 
  If $i>2$, then the sequence $W$,  $\tilde W$, $W_i$, $\cdots$, $W_n=W'$  is a chain
 from $W$ to $W'$ with length less than $n$, contradicting  to the choice of $n$.  
 If $i=1$, then $W, W_1, \tilde W$ is a $3$-cycle, contradicting to Proposition \ref{no4cycle}.
 Therefore   $i=2$ and $\tilde W\cap W_2\not=\phi$.  If $\tilde W=W_2$, then 
the sequence $W$,  $\tilde W$, $W_3$, $\cdots$, $W_n=W'$  is a chain
 from $W$ to $W'$ with length  $n-1$,  again contradicting  to the choice of $n$. 
  Thus $i=2$ and $\tilde W\not=W_2$.
But then 
 $W$, $W_1$, $W_2$, $\tilde W$  is  a  $4$-cycle, contradicting to Proposition \ref{no4cycle}.

Now $x_1\in \ol  C_v$ and $x_1\in W_0\cap W_1$ imply $W_1\cap \ol C_v\not=\phi$.
An argument similar to the one in the preceding paragraph shows  $x_2$ and 
$x_{n+1}$ lie on  the same side  of $W$.  It follows that 
  $x_2$ and $C_v$ lie on the same side  of $W$.
Therefore  an initial open segment of $x_1x_2$ lies in $\ol C_v$ and 
 does not intersect any wall
other than $W_1$. 
Choose a point $x_1'$ belonging to this  initial   open segment of  $x_1x_2$.
Now consider the sequence $W_1$,$W_2$, $\cdots$, $W_n$ and let $\sigma'$ be 
the part of $\sigma$ from $x_1'$ to 
$x_{n+1}$.  
The induction  hypothesis implies that   $\sigma'\subset \ol  C_v$.  Now the lemma follows.
    
\end{proof}

   By Theorem \ref{subm}
   $\ol  C_v/{Stab_{\Gamma}\ol  C_v}$ is compact. Lemma  \ref{connected}  implies 
$\Lambda$ is  a boundary   component of    $\ol C_v$.
   It follows that
    $\Lambda/{Stab_{\Gamma}\Lambda}$ is   a   closed $3$-manifold,  being 
      a quotient of a boundary  component of 
$\ol  C_v/{Stab_{\Gamma}\ol  C_v}$.

Next we construct a graph $G$    associated   to    $\Lambda$.    The vertex set of $G$ is in one-to-one
   correspondence with $\{W\cap \ol  C_v: W\cap \ol  C_v\not=\phi,W\in {\mathcal W^*}\}$.
   Let $v_1$ and $v_2$ be two vertices of $G$ corresponding to 
   $W_1\cap \ol  C_v$ and $W_2\cap \ol  C_v$ respectively.  There is an edge 
 connecting $v_1$ and $v_2$ if and only if $W_1\cap W_2\cap \ol  C_v\not=\phi$.

The following lemma is clear from the definitions.

\begin{Le}\label{tree}
{  The graph $G$ is a tree if there is no cycle in $X$.}
\end{Le}

\begin{Le}\label{sconnected}
{$\Lambda$ is simply connected if there is no cycle in $X$.}
\end{Le}
\begin{proof}
  We notice 
$G$ is the  nerve of the covering $\{W\cap \ol  C_v\}$ of $\Lambda$
  since $W_1\cap W_2\cap W_3=\phi$    if  $W_1, W_2, W_3$   are distinct. 
All $W\cap \ol  C_v$ and their nonempty intersections are convex and thus 
 contractible. Therefore $G$ is homotopy equivalent  to $\Lambda$. 
The lemma now follows from Lemma \ref{tree}.

\end{proof}

  Suppose  there is no cycle in $X$.  
Set $\Gamma_1=Stab_{\Gamma}\Lambda$.    Lemma \ref{sconnected}
   implies   that    $\Gamma_1\cong \pi_1(\Lambda/\Gamma_1)$    and  so 
$H\subset \Gamma_1 $  is a subgroup of the fundamental group of a
   closed  $3$-manifold.
  We may assume $\Lambda/\Gamma_1$  is orientable by replacing $\Gamma_1$ with 
  an index two subgroup
  if necessary. 
  $\Lambda/\Gamma_1$  is a Haken manifold: 
Pick any $W_1\not=W_2\in   {\mathcal W^*} $ with  $W_1\cap W_2\not=\phi$  and $W_i\cap \ol C_v\not=\phi$ ($i=1,2$);
 then $F:=W_1\cap W_2$ is a $2$-flat contained in $\Lambda$  and 
$F/Stab_{\Gamma_1}F$ is a  torus or Klein bottle    embedded  in   
$\Lambda/\Gamma_1$; the inclusion of $F/Stab_{\Gamma_1}F$  into 
$\Lambda/\Gamma_1$  clearly induces an  injective homomorphism on the fundamental groups. 
Now Theorem   \ref{haken}  implies $H$ 
   has the TA-property  since 
  any solvable subgroup of a group acting  properly and cocompactly on 
a Hadamard space is virtually free abelian.  
 The proof of Theorem \ref{maint} is now complete.

\section{Cycles of Higher Rank Submanifolds} \label{cyclesof}

Throughout this section   let 
$M=X/\Gamma$      be   a rank 1 closed real analytic 4-manifold of nonpositive sectional curvature,
  and $H\subset \Gamma$  a subgroup as described at the beginning of Section \ref{fixedpoint}.
   Recall $H$ contains no Type D isometries. 
   We proved in Section \ref{fixedpoint} that  $H$ has the TA-property if there is no cycle in $X$.   In this section 
   we discuss the Tits alternative  without assuming the     
nonexistence of cycles. 
 
   Recall   a {singular geodesic} in $X$ is a geodesic 
of the form  $\{q\}\times R\subset Q\times R=W$,   
  where   $q\in Q$   and   $W$ is a wall.  
When   $W_1$ and $W_2$ are two walls and 
$F=W_1\cap W_2$   is a $2$-flat, there are two families of parallel singular geodesics 
 in  $F$.   The angle $\alpha_F$ ($0<\alpha_F\le \frac{\pi}{2}$)  between them 
is a \emph{singular angle}.    Since modulo   
$\Gamma$ 
there are only a finite number of  walls   
  in $X$,  there are only a finite number of singular angles.

 Now we are ready to state the main 
  result of this section.

\begin{Th}\label{sameangle}
{ Let $M=X/\Gamma$    be    a rank 1 closed real analytic 4-manifold of nonpositive sectional curvature.
   Suppose there is a number $\alpha$,  $\frac{2\pi}{5}<\alpha\le \frac{\pi}{2}$  so that all the singular angles
  of $X$ are equal to $\alpha$.   If a subgroup of the fundamental group  of M is not virtually free abelian, then it contains a free 
group  of rank two.}

\end{Th}

\begin{remark}\label{example}
\emph{U. Abresch and V. Schroeder \emph{(\cite{AS})} constructed a class of real analytic 
 $4$-manifolds of nonpositive  sectional curvature,  where  all 
    the singular angles  are   $\frac{\pi}{2}$. }
\end{remark}

\subsection{Incidence Graph} \label{inci}


We  use the notation of    Section \ref{fixedpoint}.   
We first  construct a graph    $G_*$ 
which reflects    the incidence relation of the members of $ {\mathcal W^*}$.
   The vertex set of $G_*$ is in one-to-one   correspondence with  $ {\mathcal W^*}$.
  Two vertices are joined by an edge if the corresponding walls have nonempty intersection.
    $G_*$   is a connected graph.    
The vertex corresponding to the wall $W$ is still denoted    by   $W$, and the edge 
joining two vertices  $W_1$, $W_2$ is denoted   by  $W_1W_2$. 
We declare each edge of $G_*$ has length 1 and let $d_*$ be the induced path
   metric on $G_*$. 
 Since by Proposition \ref{no4cycle} there are no $n$-cycles  in $X$ 
for $n\le 4$, we have:

\begin{Le}\label{atleast5}
{  
  Any 
  injective loop in $G_*$ has length at least 5.  }
\end{Le}

The lemma in particular implies for $W_1, W_2\in  {\mathcal W^*}$ with $d_*(W_1,W_2)=2$,
   there is  a unique $W\in  {\mathcal W^*}$  with $d_*(W_1, W)=d_*(W,W_2)=1$.

    Note  
 $h({\mathcal W^*})={\mathcal W^*}$ for each   $h\in H$.  Hence   $H$ induces an action 
   on $G_*$.   
  For  $h\in H$ we denote   by  $h_*$ the isomorphism of
   $G_*$   induced by $h$.  


\begin{Le}\label{Hfixed}
{If $H$ has a fixed point in $G_*$, then $H$ has the TA-property.}
\end{Le}
\begin{proof}
Since $G_*$ is a graph, after passing to an index two subgroup if necessary,  we   may assume
  $H$ fixes a vertex $W$ of $G_*$.   Then $h(W)=W$ for all $h\in H$, here $W\subset X$ 
    denotes  the wall in  $X$.   Now the lemma follows from the proof of Lemma \ref{v}.

\end{proof}

\subsection{Tits Geodesics} \label{tits4m}

We continue  to  use the notation of    Section \ref{fixedpoint}. 
  In this section we   take a close look at  Tits   geodesics     in  
${\mathcal C}=\cup_{W\in {\mathcal W^*}}\p_T W$.

Let $W$
 be a   wall. 
A    singular  geodesic 
$c=\{q\}\times R\subset Q\times R=W$ determines two points   
$c(+\infty)$, $c(-\infty)$ in $\partial_\infty W$.   Set $w(+\infty)=c(+\infty)$,
$w(-\infty)=c(-\infty)$   
  and call $w(+\infty)$, 
$w(-\infty)$    the \emph{poles}   of  $W$.

The following proposition follows from  the results in \cite{HS1}.

\begin{Prop}\label{w1w2}
{Let $W_1, W_2\in {\mathcal W^*}$ be two walls.  \newline
\emph{(i)} If  $d_*(W_1, W_2)=1$,  i.e., if   $W_1\cap W_2=F$ is a  $2$-flat,    then 
 $\partial_\infty W_1\cap \partial_\infty W_2=\partial_\infty F$;\newline
\emph{(ii)} If $d_*(W_1, W_2)=2$,      then 
$\partial_\infty W_1\cap \partial_\infty W_2= \{w(+\infty), w(-\infty)\}$
    where $W$ is the unique wall with $d_*(W_1, W)=d_*(W,W_2)=1$;\newline
\emph{(iii)} If $d_*(W_1, W_2)\ge 3$,      then 
$\partial_\infty W_1\cap \partial_\infty W_2= \phi$.}
\end{Prop}




Let $W_1, W_2\in {\mathcal W^*}$ with $d_*(W_1, W_2)=1$.  Then 
   $W_1\cap W_2=F$  is  a $2$-flat  and $\partial_T W_1\cap \partial_T W_2=\partial_T F$ is 
  isometric to the unit circle.  
$\partial_T F$    admits  a   unique metric graph structure   with   vertex set  
$\{w_1(+\infty),   w_1(-\infty),   w_2(+\infty),    w_2(-\infty)\}$.  

 Let $W=Q\times R\in   {\mathcal W^*}$.    Since  $W$ is closed and $Q$  is a nonflat
 Hadamard $2$-manifold,     
$\partial_T W$   
     admits   a unique metric graph structure   with the following  properties:\newline 
(a) the vertex set  is $\{w'(+\infty), w'(-\infty):  d_*(W', W)\le 1\}$;\newline
(b) for each wall $W'$ with $d_*(W', W)=1$, the inclusion
   $\partial_T(W'\cap W)\subset \partial_T W$   is an isometric embedding 
 between  metric graphs;\newline
(c) each edge connecting the two poles $w(+\infty)$,  $w(+\infty)$ has length $\pi$.

   It  follows from  the results in \cite{HS1}  that 
$\mathcal C$   admits  a  unique metric graph structure
    with   vertex set  $\{w(+\infty), w(-\infty):  W\in  {\mathcal W^*}\}$
     such that for each $W\in  {\mathcal W^*}$, the inclusion  $\p_T W\subset {\mathcal C}$ is an isometric embedding 
  between metric graphs.

We next  look at how Tits geodesics  in $\mathcal C$ travel between  different 
   $\partial_\infty W$,
$W\in {\mathcal W}^*$.

\begin{Le}\label{noback}
{Let $W\in {\mathcal W}^*$, $\sigma: [a, b]\rightarrow \partial_T X$ a minimal 
       geodesic
and $t_0\in (a, b)$.   If  there is some   $\epsilon>0$ such that 
$\sigma(t)\in \partial_\infty W$ for  $t\in (t_0-\epsilon, t_0]$
 and $\sigma(t)\notin \partial_\infty W$ for  $t\in (t_0, t_0+\epsilon)$,
then $\sigma(t)\notin \partial_\infty W$ for $t\in (t_0, b]$.}
\end{Le}
\begin{proof}
Suppose the lemma is false and let $t_1=\min\{t\in (t_0, b]: \sigma(t)\in \partial_\infty W\}$.  
  Notice the minimal in the definition of $t_1$ makes sense since $\partial_\infty W$ is closed
  in   $\partial_\infty X$.
     By definition  $\partial_\infty W\cap \sigma_{|(t_0,t_1)}=\phi$.    
     Since  $\sigma(t_0),\sigma(t_1)\in \partial_\infty W$ we have $d_T(\sigma(t_0),\sigma(t_1))\le \pi$.
  Let $\sigma'\subset \partial_T W$ be a minimal geodesic from $\sigma(t_0)$ to $\sigma(t_1)$. 
  Then $\sigma'\cup \sigma_{|[t_0,t_1]}$ is a closed geodesic in the  $CAT(1)$ space
  $\partial_T X$ with length 
$$length(\sigma')+length(\sigma_{|[t_0,t_1]})=2 d_T(\sigma(t_0),\sigma(t_1))\le 2\pi.$$
  It follows that $d_T(\sigma(t_0),\sigma(t_1))=\pi$.  Since    
$\sigma_{|(t_0-\epsilon, t_0]}\subset \partial_\infty W$,    the description of  $\p_T W$ shows 
$d_T(\sigma(t), \sigma(t_1))<\pi$ for $t\in (t_0-\epsilon, t_0)$, contradicting to 
the fact that $\sigma$ is minimal.

\end{proof}

 The  following proposition   follows from   Lemma \ref{noback}  and  the description  of ${\mathcal  C}$.

\begin{Prop}\label{inj}
{Given any minimal geodesic $\sigma:[a,b]\rightarrow {\mathcal  C}$, 
  there  are  walls   $W_1$, $W_2$, $\cdots$, $W_n $ and  numbers
  $a< t_1<t_2\cdots < t_n< b$ with the following properties:\newline
\emph{(i)} $\{\sigma(t_i), 1\le i\le n\}$  is the set of  poles  in the interior  of   $\sigma$;\newline
\emph{(ii)} $\sigma(t_i)$ is a pole of $W_i$ for each $1\le i\le n$;\newline
\emph{(iii)}  the sequence of walls $W_1, W_2, \cdots, W_n $   determines 
   an injective 
edge path in $G_*$.}
\end{Prop}

\subsection{Subgroups Containing Type C Isometries} \label{gctypec}

In this section we     show 
     the subgroup $H\subset \Gamma$  has the TA-property if the singular  angles 
 are large and  $H$ contains  a Type C isometry.

 We recall  there are only  a finite number of singular angles.

\begin{Prop}\label{typecgroup}
{ Let $M=X/\Gamma$    be    a rank 1 closed real analytic 4-manifold of nonpositive sectional curvature.
   Suppose  all the singular angles
  of $X$ are   strictly larger than $\frac{\pi}{3}$.  
 If $H\subset \Gamma $ is  a subgroup  as described  at the beginning of Section \ref{fixedpoint}
    and contains a Type   C isometry,  then  $H$ 
 has the TA-property.}
\end{Prop}

\begin{proof}
    For any two Type C isometries $h, \tilde h\in H$,
  we set $n(h, \tilde h)=d_*(W, \tilde W)$ where $W=P(h)$ and $\tilde W=P(\tilde h)$. 
Notice $h(+\infty),  h(- \infty)$ are the poles of $W$, and 
 $\tilde h(+\infty)$,  $\tilde h(- \infty)$ are the poles of $\tilde W$. 
First suppose there are two Type C isometries $h$ and $\tilde h$ with 
$n(h,\tilde h)\ge 3$.    
 For any $\xi\in \{ h(+\infty),  h(- \infty)\}$ and $\eta\in \{ \tilde h (+\infty),  \tilde h(- \infty)\}$,  let
  $\sigma: [0,a]\rightarrow \partial_T X$ be a minimal     geodesic from $\xi$ to $\eta$.
If there is no pole 
 in the interior of $\sigma$,  then $\sigma\subset \partial_\infty W\cap  \partial_\infty {\tilde W}$, 
contradicting to Proposition \ref{w1w2}. 
Let $\sigma(t_i)$ ($i=1,2,\cdots, n$)  with  $0<t_1<t_2\cdots < t_n<a$
be all the poles in the interior of $\sigma$. 
  If $n=1$, then $\sigma(t_1)\in  \partial_\infty W\cap  \partial_\infty {\tilde W}$,
  again contradicting to Proposition \ref{w1w2}. Therefore there are at least
 two poles in the interior of $\sigma$.  Since the endpoints 
of   $\sigma$ are also poles,  $\sigma$    contains at least 4 poles. By assumption 
on the singular angles,  the Tits distance between any two   distinct  
poles is   $>\frac{\pi}{3}$.
  It follows that $d_T(\xi, \eta)>\pi$.  By  Theorem \ref{kim}  
  $<h,\tilde h>$   contains  a free group of rank two.

Now suppose there are two Type C isometries $h$ and $\tilde h$ with 
$n(h,\tilde h)=2$.  Let $W_1=P(h)$, $W_2=P(\tilde h)$.  
  There exists an integer $k\ge 1$  with 
$Min(h^k)=W_1$ and $Min({\tilde h}^k)=W_2$.  
 By Lemma \ref{atleast5}
  there is a unique  wall 
  $W$ with $W\cap W_1\not=\phi$  and $W\cap W_2\not=\phi$.   We have
 $h^k(W)=W$, ${\tilde h}^k(W)=W$. Let $W\cap W_1=c_1\times R\subset Q\times R=W$ and
 $W\cap W_2=c_2\times R\subset Q\times R=W$,  where  $c_1, c_2$ are complete geodesics in $Q$.
   Theorem  \ref{subm}
  implies  $\{c_1(+\infty), c_1(-\infty)\}\cap \{c_2(+\infty), c_2(-\infty)\}=\phi$.
 $h^k$ acts on $W$ as $h^k=(h_1, t): Q\times R\rightarrow Q\times R$,  where 
  $h_1:Q\rightarrow Q$ is a hyperbolic isometry of $Q$ with $c_1$ as an axis and $t$ is a translation
  of $R$. 
 Similarly  ${\tilde h}^k$   acts on $W$ as ${\tilde h}^k=(h_2, t')$  where $h_2$ is 
a hyperbolic isometry of $Q$ with $c_2$ as an axis  and $t'$ is a translation of $R$.
  Since $Q$ is nonflat and  admits a cocompact group of isometries, $Q$ is hyperbolic in the sense of Gromov. 
   Theorem  \ref{dyn}  implies  $<h_1, h_2>$  contains a  free  group of rank two. 
   It follows that 
     $<h^k,{\tilde h}^k>$  contains a  free  group of rank two.

Now suppose  $n(h,\tilde h)\le 1$  for any two Type C isometries $h, \tilde h\in H$.
For any three Type C isometries $h_1,h_2, h_3\in H$, let $W_i=P(h_i)$  ($i=1,2,3$).
   If $W_1\not=W_2$ and $W_1\not=W_3$,   then $W_2=W_3$ since   $n(h_2,h_3)\le 1$
  and   $G_*$  has  no injective edge loop with length 3.   It  follows  that 
  the set $\{ P(h):h\in H \text{ is of Type C}\}$ consists of one wall or two
 intersecting walls. Notice the set $\{ P(h):h\in H \text{ is of  Type C}\}$ 
 is invariant  under the action of $H$ since the conjugate of a Type C isometry
is still a Type C isometry and $P(ghg^{-1})=g(P(h))$.
 The proposition  now follows from Lemma \ref{Hfixed}.

\end{proof}

\subsection{Admissible Subsets} \label{admissible}

For a  Hadamard space $X$, 
    $\xi\in \partial_T X$ and $r>0$, we let  $\ol{B}(\xi, r)$  
 be the closed metric ball with center
  $\xi$ and radius $r$:  
$\ol{B}(\xi, r)=\{\eta\in \partial_T X:   d_T(\xi,\eta)\le r\}$. 

\begin{Def}\label{adm}
{Let  $X$  be  a Hadamard space and  $H\subset  Isom(X)$. 
   A    nonempty   
$H$-invariant   subset $M\subset \partial_\infty X$
  is    \emph{$H$-admissible} if 
 $M\subset \ol{B}(h(+\infty), \pi)\cap \ol{B}(h(-\infty), \pi)$ for  each 
   hyperbolic isometry   $h\in H$.}
\end{Def}

\begin{Prop}\label{empadm}
{Let 
$M=X/\Gamma$      be   a rank 1 closed real analytic 4-manifold of nonpositive sectional curvature,
  and $H\subset \Gamma$  a subgroup as described at the beginning of Section \ref{fixedpoint}. 
  Then there exists a  $H$-admissible subset.}
\end{Prop}
\begin{proof}
We notice  $h(+\infty),  h(-\infty)   \in  \Lambda(H)$
for any $h\in H$.   Let $M  \subset \Lambda(H)$ be a    
 $H$-minimal set.  
   Since by assumption  $H$ does not contain any rank one  isometries,
    Proposition \ref{dis} implies    $M$ is    $H$-admissible.   

\end{proof}

\subsection{Intersection of Tits Balls} \label{titsballs}

In order to study $H$-admissible subsets in   $\partial_\infty X$,
   we need to look at the intersection of  
$\ol{B}(h(+\infty),\pi)$ and $\ol{B}(h(-\infty),\pi)$
  for  $h\in H$.

Call $h\in H$   a  \emph{squared Type B isometry} if $h=g^2$ for a Type B isometry $g\in H$.
 Clearly  $h\in H$ is a  squared Type B isometry  if and only if any of its conjugates is 
a  squared Type B isometry.  
 Notice for a  squared Type B isometry  $h$,  $Min(h)=P(h)$. 
If  $h\in H$ is of  Type A, 
   then $Fix(h_*)$ is the vertex $W$ in $ G_*$,
 where $W$ is the only wall containing $Min(h)\subset X$. If 
$h\in H$ is a  squared Type B isometry,   then $Fix(h_*)$ is the edge $W_1W_2$ in  $G_*$,
where $W_1$, $W_2$ are the two walls with $W_1\cap W_2=Min(h)\subset X$.
For any $h\in H$ and any wall $W$, set 
$$\ol{B}(h,W)=\ol{B}(h(+\infty),\pi)\cap \ol{B}(h(-\infty),\pi)\cap \partial_\infty W.$$

\begin{Le}\label{inter3}
{Suppose the singular angles of $X$ are all  equal  to  $\alpha>\frac{2}{5}\pi$. 
Let $h\in H$  be a Type A or squared Type B isometry,  and 
 $W\in {\mathcal W}^*$ with $d_*(W, Fix(h_*))\ge 3$.  Then
  $ \ol{B}(h,W)   =\phi$.}
\end{Le}
\begin{proof}
Suppose the lemma is false and pick $\xi\in \ol{B}(h,W)$.
Let $c:[0,a] \rightarrow \partial_T X$ ($a\le \pi$) be a  minimal     geodesic from 
$ h(+\infty)$ to $\xi$, $W_1, \cdots, W_k$  a sequence of walls and  $\{t_i\}$  ($0<t_1<\cdots<t_k<a$)
 a sequence of numbers   as in  Proposition \ref{inj}.
By definition of Type A and squared Type B isometries,
 $h(+\infty)$ is not a pole and $h_*(W_1)=W_1$.  $a\le \pi$ implies 
$k\le 3$ since by assumption the Tits distance between any two   distinct  poles
is at least $\frac{2}{5}\pi$.  If  there is no pole in the interior of $c$, 
  then there is a wall $\tilde W$ with 
   $h_*(\tilde W)=\tilde W$  and $d_*(\tilde W, W)\le 1$,   
  contradicting to the assumption 
  $d_*(W, Fix(h_*))\ge 3$.  If $k=1$, then $\xi\in \partial_\infty {W_1}\cap \partial_\infty W$
  and by Proposition \ref{w1w2} we have $d_*(W_1, W)\le 2$, again contradicting to
the assumption.  If $k=2$, then $\xi\in \partial_\infty {W_2}\cap \partial_\infty W$
    and since  $d_*(W, Fix(h_*))\ge 3$ 
  we have  $d_*(W_2, W)=2$ and   $\xi$ must be a pole. In this case let 
$W_3$ be the wall with $\xi$ as one of its two poles and set $t_3=a$.
If $k=3$, then $\xi$ can not be a pole since $a\le \pi$ and the Tits distance between any two poles
is at least $\frac{2}{5}\pi$.
   In any case  we have  $\xi\in B(c(t_3), \frac{\pi}{5})$.   Similarly
  if  $c':[0,b] \rightarrow \partial_T X$ ($b\le \pi$) is  a  minimal     geodesic from 
$ h(-\infty)$ to $\xi$,
 we have  walls $W_1', W_2', W_3'$ ($h_*(W_1')=W_1'$)   and poles  $c'(t_1')$, $c'(t_2')$, $c'(t_3')$
     such that  $\xi\in  B(c'(t_3'), \frac{\pi}{5})$.  
It follows from the  triangle inequality that $d_T(c(t_3),c'(t_3'))<\frac{2}{5}\pi$. 
  The assumption on singular angles implies $c(t_3)=c'(t_3')$.  
Consequently
  $W_3'=W_3$.  

We next show  $W_1'=W_1$. Note  $t_1, t_1'<\frac{\pi}{5}$.   Since  $d_T(h(+\infty),h(-\infty))=\pi$,  
   triangle  inequality   implies  $d_T(c(t_1), c'(t_1'))>\frac{3}{5}\pi$.     The two poles 
$c(t_1)$, $c'(t_1')$  both lie on  $\partial_\infty {W_1}$.
   Notice  for any pole $\eta$  in  $\partial_\infty {W_1}$,  we have  
  $d_T(\eta, c(t_1))\in \{\alpha, \pi-\alpha, \pi\}$. 
The assumption on $\alpha$ now implies $d_T(c(t_1), c'(t_1'))=\pi$.  Therefore $c'(t_1')$ 
  is also a pole of $W_1$, and $W_1'=W_1$.

   Now we have two injective edge paths of length 2 from $W_1=W_1'$ to
$W_3=W_3'$:  $W_1W_2W_3$, $W_1'W_2'W_3'$.  By Lemma \ref{atleast5} we  have $W_2=W_2'$.
   Notice $t_3-t_2=\alpha$ otherwise $t_3-t_2\ge \pi-\alpha$ and 
$a>(t_3-t_2)+(t_2-t_1)\ge (\pi-\alpha)+\alpha=\pi$.   The same argument also shows 
$\alpha<\frac{\pi}{2}$ and $t_3'-t_2'=\alpha$.  
 Now the three poles $c(t_2)$, $c'(t_2')$, $c(t_3)$ all lie on  the circle
 $\partial_\infty {W_2}\cap \partial_\infty {W_3}$ and  $d_T(c(t_2), c(t_3))=d_T(c'(t_2'), c(t_3))=\alpha$.
Since $\alpha<\frac{\pi}{2}$  we have $c'(t_2')=c(t_2)$.   Similarly we conclude
$c'(t_1')=c(t_1)$.   

Since  $t_1, t_1'<\frac{\pi}{5}$, by triangle 
 inequality 
      $$d_T(h(+\infty),h(-\infty))\le d_T(h(+\infty), c(t_1))+d_T(c'(t_1'),h(-\infty))<\frac{2}{5}\pi,$$
   contradicting to the fact $d_T(h(+\infty),h(-\infty))=\pi$.  

\end{proof}

    If   $h\in H$ is of  Type A  and $W$ is a wall with  $d_*(W, Fix(h_*))=2$, there is exactly
 one wall $W'$ with $d_*(W, W')=d_*(W', Fix(h_*))=1$.  
If $h\in H$ is a squared  Type B   isometry  and  $W$ is a wall with  $d_*(W, Fix(h_*))=2$, 
  then the set 
$$\{W': d_*(W, W')=d_*(W', Fix(h_*))=1\}$$  
consists of either one or two elements.  These assertions follow from Lemma \ref{atleast5}.

\begin{Le}\label{inter2}
{ Let $h\in H$  be a Type A or  squared Type B isometry, and 
 $W\in {\mathcal W}^*$   such that  the following   holds: $d_*(W, Fix(h_*))=2$. 
Suppose the singular angles of $X$ are all  equal  to  $\alpha>\frac{2}{5}\pi$.  
\newline
\emph{(i)} If  there is only    one wall
$W'$ with  
 $d_*(W, W')=d_*(W', Fix(h_*))=1$, then 
$$\ol{B}(h,W)\subset \partial_\infty W  \cap \partial_\infty {W'}-\{w(+\infty), 
      w(-\infty)\};     $$
\emph{(ii)} If  there  
are two  walls $W_1$, $W_2$ 
 with 
$d_*(W, W_i)=d_*(W_i, Fix(h_*))=1,$
     then 
  $$ \ol{B}(h,W)        \subset 
(\partial_\infty W\cap \partial_\infty {W_1})\cup (\partial_\infty W\cap \partial_\infty {W_2})       -\{w(+\infty), 
      w(-\infty)\}.     $$}
\end{Le}
\begin{proof}
  We fix an arbitrary 
$\xi  \in  \ol{B}(h,W)$.
Let $c:[0,a] \rightarrow \partial_T X$ be a  minimal     geodesic from 
$ h(+\infty)$ to $\xi$, and $c':[0,b] \rightarrow \partial_T X$  a  minimal     geodesic from 
$ h(-\infty)$ to $\xi$, where $a, b\le \pi$.  
  Let 
$W_1, \cdots, W_k$ be a sequence of walls and  
  $\{t_i\}$ ($0<t_1<\cdots<t_k<a$) 
  a sequence of numbers    provided  by  Proposition \ref{inj}
   corresponding to $c$, and  $W_1', \cdots, W'_{k'}$
 and  $\{t'_i\}$  ($0<t'_1<\cdots<t'_{k'}<b$)  
 corresponding to $c'$. As in the proof of Lemma \ref{inter3} we see 
$h_*(W_1)=W_1$, $h_*(W'_1)=W'_1$  and $1\le k\le 3$,  $1\le  k'\le 3$.

   If $k=1$, then $\xi\in \partial_\infty {W_1}\cap \partial_\infty W$. Since $d_*(W, Fix(h_*))=2$ and $h_*(W_1)=W_1$, 
 Proposition \ref{w1w2}  implies $\xi$ is a pole of $W'$, where $W'$ is the only wall 
with $d_*(W, W')=1$ and  $d_*(W', W_1)=1$.  We clearly have $\xi\in \partial_\infty W\cap \partial_\infty {W'}-\{w(+\infty), 
      w(-\infty)\}$ in this case.  From now on we assume $2\le k,k'\le 3$. 

 Suppose $k=2$.     In this case  $\xi\in \partial_\infty {W_2}\cap \partial_\infty W$. If $\xi$ is not a pole, 
then $W_2$ is the only wall with  $d_*(W, W_2)=d_*(W_2, W_1)=1$ and clearly
  $$\xi\in \partial_\infty W\cap \partial_\infty {W_2}-\{w(+\infty), w(-\infty)\}.$$
  Suppose $\xi$ is   a  pole of   some   wall $W_3$. Then the 
 assumption on $\alpha $ implies $k'=2$.  In this case we have two injective
 edge paths in $G_*$: $W_1W_2W_3$ and $W'_1W'_2W_3$.  The proof of Lemma \ref{inter3}
 yields a contradiction if $W'_1=W_1$.     Hence   $W'_1\not=W_1$,    
$h$   is a Type B isometry  and    $Fix(h_*)=W_1W'_1$.
 Notice $t_1<\frac{\pi}{5}$
   since $a\le \pi$ and $c(t_1)$, $c(t_2)$, $\xi$ are three poles on $c$.  Similarly
  $t'_1<\frac{\pi}{5}$.
    Recall $c(t_1)$ is a pole of $W_1$. Let $p\not=c(t_1)$ be the other pole of $W_1$.
  Since $h(+\infty)$, $h(-\infty)$, $p$ and $c(t_1)$  all lie on the circle 
  $\p_T W_1\cap \p_T W_1'$,  we see 
     $d_T(h(-\infty), p)=d_T(h(+\infty), c(t_1))=t_1$.  It follows that
  $$d_T(c'(t'_1), p)\le d_T(c'(t'_1),h(-\infty))+d_T(h(-\infty), p)=t'_1+t_1<\frac{2}{5}\pi.$$
    As  $c'(t'_1)$  and  $p$ are poles of two distinct walls $W'_1$  and  $W_1$  respectively, 
     we  have  
   $ c'(t'_1)\not=p$
 and  $d_T(c'(t'_1), p)\ge \alpha> \frac{2}{5}\pi$, a contradiction. 

   Now we assume $k=k'=3$.   
In this case $\xi$ is not a pole. The proof of Lemma \ref{inter3} shows $c(t_3)=c'(t'_3)$ and 
 $W_3=W'_3$. If $W'_1=W_1$ then the proof of Lemma \ref{inter3}
  yields a contradiction.   So $W'_1\not=W_1$,   $Fix(h_*)=W_1W'_1$    and 
   $h$ is  a Type B isometry.  
  Now consider $c(t_3)=c'(t'_3)\in \ol{B}(h,W_3)$ instead of $\xi$ and the preceding 
paragraph
 yields a contradiction.

\end{proof}

\subsection{Proof of Theorem \ref{sameangle} }\label{proofof}

    In this section we finish the proof of  Theorem \ref{sameangle}. 
  By Lemma \ref{Hfixed} and Proposition \ref{typecgroup}, we may assume 
$H$  does not have a global fixed point in $G_*$ and 
each  nontrivial $h\in H$ is of Type A or Type B.  

Let $M\subset \partial_\infty X$ be a  $H$-admissible subset. 
By Lemma \ref{inter3},
  if  $W$ is a wall    with  $d_*(W, Fix(h_*))\ge 3$ for some $h\in H$, then 
$M\cap \partial_\infty W=\phi$. 

\begin{Le}\label{a0}
{Let $M\subset \partial_\infty X$ be     $H$-admissible.  
Suppose   $W$  is  a wall   such that   $Fix(h_*)=\{W\}$
for  a Type A isometry
$h\in H$. 
     Then $M\cap \partial_\infty W=\phi$.}

\end{Le}
\begin{proof}
Since $Fix((ghg^{-1})_*)=g_*(Fix(h_*))$ and 
the conjugate of a Type A isometry is  still   a Type A isometry,
   the assumption that $H$ does not have a global fixed point
 in $G_*$ implies that there is a wall $W'\not=W$ with $Fix(h'_*)=\{W'\}$ for a Type A isometry
  $h'\in H$.   If $d_*(W', W)\ge 3$, then Lemma \ref{inter3} implies 
$\ol{B}(h',W)=\phi$.  By the definition of a  $H$-admissible set, 
$M\subset \ol{B}(h'(+\infty),\pi)\cap \ol{B}(h'(-\infty),\pi)$
  and thus
 $M\cap \partial_\infty W\subset \ol{B}(h'(+\infty),\pi)\cap \ol{B}(h'(-\infty),\pi)\cap \partial_\infty W=\ol{B}(h',W)=\phi$.
From now on we assume    $d_*(W', W)\le 2$.

Suppose $d_*(W', W)=1$.  Then $h'(W)\not=W$ and $d_*(W, W')=d_*(h'(W), W')=1$.  
   By  Lemma \ref{atleast5}  $d_*(W, h'(W))=2$.  Thus we may assume $d_*(W', W)=2$.    Then 
$d_*(h(W'), W)=2$.  Let $W_1$, $W_2$ be the  unique walls with 
  $d_*(W, W_1)=d_*(W_1, W')=1$,     $d_*(W, W_2)=d_*(W_2, h(W'))=1$.  
Notice $W_2=h(W_1)\not=W_1$. 
  Lemma \ref{inter2}
 implies  
$$\ol{B}(h', W)\subset \partial_\infty W\cap \partial_\infty {W_1}-\{w(+\infty), w(-\infty)\}$$
  and 
$$\ol{B}(hh'h^{-1}, W)\subset \partial_\infty W\cap \partial_\infty {W_2}-\{w(+\infty), w(-\infty)\}.$$
   $M$ is $H$-admissible implies

$\;\;M\cap \partial_\infty W\subset   \ol{B}(h', W)\cap \ol{B}(hh'h^{-1}, W)\subset $ 

$\;\;\subset (\partial_\infty W\cap \partial_\infty {W_1}-\{w(+\infty), w(-\infty)\})\cap (\partial_\infty W\cap \partial_\infty {W_2}-\{w(+\infty), w(-\infty)\})\subset $

$\;\;\subset \partial_\infty {W_1}\cap \partial_\infty {W_2}-\{w(+\infty), w(-\infty)\}=\phi.$

\noindent
 The  last  equality follows from 
 Proposition \ref{w1w2} and the facts that  $W_1\not=W_2$ and $d_*(W_1, W)=1$,  $d_*(W, W_2)=1$.   

\end{proof}

\begin{Le}\label{a1}
{Let $M\subset \partial_\infty X$ be    $H$-admissible.  Suppose  
$W$  is a  wall   with   $d_*(Fix(h_*),W)=1$  for a Type A isometry 
 $h\in H$.  
      Then $M\cap \partial_\infty W=\phi$.}
\end{Le}
\begin{proof}
Suppose $M\cap \partial_\infty W\not=\phi$.   
Let $Fix(h_*)=\{W_0\}$. The argument in the  proof of Lemma \ref{a0}  shows that there is a wall
$W'$  and a Type A isometry  $h'\in H$ such that
 $Fix(h'_*)=\{W'\}$ 
 and  $d_*(W_0, W')=2$ or $3$.  

Let  us first assume $d_*(W_0, W')=3$.
 Then $M\cap \partial_\infty W\not=\phi$  implies $d_*(W, W')=2$  and   $d_*(W, h(W'))=2$.
Let $W_1$ be the unique wall with $d_*(W, W_1)=d_*(W_1, W')=1$,
 and $W_2$ the unique wall with $d_*(W, W_2)=d_*(W_2, h(W'))=1$. 
 Lemma \ref{inter2} implies

$\;\;\phi\not= M\cap \partial_\infty W\subset\ol{B}(h', W)\cap \ol{B}(hh'h^{-1}, W)\subset$

$\;\;\subset (\partial_\infty W\cap \partial_\infty {W_1}-\{w(+\infty), w(-\infty)\})\cap (\partial_\infty W\cap \partial_\infty {W_2}-\{w(+\infty), w(-\infty)\})\subset $

$\;\;\subset \partial_\infty {W_1}\cap \partial_\infty {W_2}-\{w(+\infty), w(-\infty)\}$, 

\noindent
and therefore $W_1=W_2$.
 Since $M$ is $H$-invariant, $M\cap \partial_\infty W\not=\phi$  implies $M\cap \partial_\infty{h(W)}$
    is nonempty.
  The above argument  applied to $h(W)$ instead of $W$ shows there exists   a wall 
$W'_1$ with $d_*(h(W), W'_1)=d_*(W'_1, W')=d_*(W'_1, h(W'))=1$.  
  If $W'_1\not=W_1$,  then $W_1W'W'_1h(W')W_1$ is an injective edge loop
with length 4, a contradiction. If $W'_1=W_1$, then $W_0WW_1h(W)W_0$ is  
an injective edge loop
with length 4, again a  contradiction.

  Now assume $d_*(W_0, W')=2$.  
$M\cap \partial_\infty W\not=\phi$, Lemma \ref{a0} and the remark before  Lemma \ref{a0}
  imply   $d_*(W, W')=1$ or $2$. 
By replacing $W'$ with $h(W')$ if necessary we 
 may assume  $d_*(W, W')=2$.    The argument in the preceding paragraph  again  yields
a contradiction.

\end{proof}

\begin{Le}\label{a2}
{Let $M\subset \partial_\infty X$ be    $H$-admissible.  Suppose  
$W$  is a  wall  with  $d_*(Fix(h_*),W)=2$  for a Type A isometry 
 $h\in H$.  
      Then $M\cap \partial_\infty W=\phi$.}

\end{Le}
\begin{proof} 
   Let $Fix(h_*)=\{W_0\}$ and $W_1$ be the unique wall  with
 $d_*(W_0, W_1)=1$ and  $d_*(W_1, W)=1$.  
 Then by Lemma \ref{inter2}
$$M\cap \partial_\infty W\subset \ol{B}(h, W)\subset \partial_\infty {W_1}\cap \partial_\infty W-\{w(+\infty), 
      w(-\infty)\}\subset \partial_\infty {W_1}.$$
  But Lemma \ref{a1} implies 
$M\cap \partial_\infty {W_1}=\phi$. It follows $M\cap \partial_\infty W\subset M\cap \partial_\infty {W_1}=\phi$.

\end{proof}

    Lemmas \ref{a0}, \ref{a1}, \ref{a2}, Proposition \ref{empadm} and the remark before 
Lemma \ref{a0} together 
imply that $H$  contains no Type A isometries. 
    From now on we  assume 
each nontrivial $h\in H$ is of Type B  and  $H$  does not have a global fixed point in $G_*$. 
Recall $h\in H$ is a  squared Type B isometry  if and only if any of its conjugates is 
a  squared Type B isometry.

\begin{Le}\label{hh'}
{If there are  squared Type B isometries $h,h'\in H$     with
$d_*(Fix(h_*), Fix(h'_*))\ge 3$, then $H$ contains  a  free group of rank two.}

\end{Le}
\begin{proof} 
Let $\xi\in \{h(+\infty), h(-\infty)\}$,  $\eta\in \{h'(+\infty), h'(-\infty)\}$
  and $\sigma:[0, a]\rightarrow \partial_T X$ be a minimal     geodesic from 
 $\xi$ to $\eta$.  Let the walls $W_1, \cdots, W_k$ and the numbers $0<t_1< \cdots< t_k<a$
   be provided by Proposition  \ref{inj}.  Then $h_*(W_1)=W_1$ and $h'_*(W_k)=W_k$.
Since $d_*(Fix(h_*), Fix(h'_*))\ge 3$,  we have $k\ge 4$.   Since each singular angle
  equals $\alpha>\frac{2}{5}\pi$, the length of $\sigma$ is at least $3\times \frac{2}{5}\pi>\pi$.
 So $d_T(\xi, \eta)>\pi$ and the lemma now follows from  Theorem \ref{kim}.

\end{proof}

\begin{Le}\label{lin}
{Suppose each nontrivial $h\in H$ is of Type B  and
 $H$  does not have a global fixed point in $G_*$.
  Then there are   squared Type B isometries  $h, h'\in H$ with
$d_*(Fix(h_*), Fix(h'_*))\ge 2$.}

\end{Le}
\begin{proof}
First suppose $d_*(Fix(h_*), Fix(h'_*))=0$  for all squared Type  B isometries $h, h'\in H$.
  Fix a  squared Type B isometry $h\in H$  and let $Fix(h_*)=W_1W_2$.    Since 
$H$  does not fix the midpoint of  $W_1W_2$,  there is  some   element  $g$  in  $ H$ such that
$g_*(W_1W_2)\cap W_1W_2=\{W_1\}$ or $\{W_2\}$.  After possibly relabeling 
 $W_1$ and $W_2$ we may assume $g_*(W_1W_2)\cap W_1W_2=\{W_1\}$. 
 Set $h_1=h$ and $h_2=ghg^{-1}$.  $h_1$ and $h_2$ are two squared Type B isometries with
 $Fix({h_1}_*)\cap Fix({h_2}_*)=\{W_1\}$.   Since 
$H$  does not fix $W_1$,  there is some $k\in H$ such that $k_*(W_1)\not=W_1$. 
 Set $k_1=kh_1k^{-1}$, $k_2=kh_2k^{-1}$.  Then  $k_1$ and $k_2$ are two squared Type B isometries with
 $Fix({k_1}_*)\cap Fix({k_2}_*)=\{k_*(W_1)\}$. 
      Now  it is not hard to derive from   Lemma \ref{atleast5}  that  
there are $i,j\in \{1,2\}$   such  that   $Fix({h_i}_*)\cap Fix({k_j}_*)=\phi$, a contradiction.

Now assume   $d_*(Fix(h_*), Fix(h'_*))\le 1$  for any two  squared Type B isometries  $h, h'\in H$.
 By the   above paragraph, there   are    two squared Type B isometries  $g, h\in H$   with 
$d_*(Fix(h_*), Fix(g_*))=1$.  Let $Fix(h_*)=W_2W_3$, 
$g_*(W_1)=W_1$ with $d_*(W_1, W_2)=1$.     Notice 
$Fix((ghg^{-1})_*)=g(W_2)g(W_3)$.  
    It follows from     $d_*(Fix(h_*), Fix((ghg^{-1})_*))\le 1$   and  Lemma \ref{atleast5}
   that  $d_*(W_3, g(W_3))=1$.
Similarly  we have  $d_*(W_3, g^2(W_3))=1$  and  $d_*(g^2(W_3), g(W_3))=1$. It follows that we have 
 an injective edge loop   $W_3g(W_3)g^2(W_3)W_3$ with length 3, a contradiction.  

\end{proof}

\begin{Le}\label{b0}
{Let  $M\subset \partial_\infty X$ be   $H$-admissible, 
  and 
 $h,g \in H$   be    two   squared Type B  isometries  with 
  $d_*(Fix(h_*),Fix(g_*))=2$  and $Fix(h_*)=WW'$.  
  Then   $M\cap \partial_\infty W= M\cap \partial_\infty {W'}   =\phi$.}
\end{Le}
\begin{proof}
We  shall prove $M\cap \partial_\infty W=\phi$. The proof of $M\cap \partial_\infty {W'}=\phi$  is similar.
If  $d_*(W,Fix(g_*))=3$,  then by Lemma \ref{inter3} $M\cap \partial_\infty W \subset \ol{B}(g,W)=\phi$.
Suppose $d_*(W,Fix(g_*))=2$.  
  Then  Lemma \ref{inter2} implies  that there is a set of walls 
$\{W_i\}_{i=1}^N$  ($N=1$ or $2$)   with  $d_*(W,W_i)=d_*(W_i, Fix(g_*))=1$   and 
$$M\cap \partial_\infty W\subset \ol B(g,W)\subset \cup_{i=1}^N (\partial_\infty W\cap \partial_\infty {W_i})-\{w(+\infty), w(-\infty)\}.$$
 Fix an integer $j$   with  $\{W_1, W_N\}\cap \{h^j(W_1), h^j(W_N)\}=\phi$. 
  Since  $ M\cap \partial_\infty W$ is  also  contained in  
$\ol B(h^jgh^{-j},W)\subset \cup_{i=1}^N (\partial_\infty W\cap \partial_\infty {h^j(W_i)})-\{w(+\infty), w(-\infty)\}$,   
it follows  
$ M\cap \partial_\infty W$  is   contained in the intersection of 
$\cup_{i=1}^N (\partial_\infty W\cap \partial_\infty {W_i})-\{w(+\infty), w(-\infty)\}$    and 
$\cup_{i=1}^N (\partial_\infty W\cap \partial_\infty {h^j(W_i)})-\{w(+\infty), w(-\infty)\}$,
   which is empty.

\end{proof}

\begin{Le}\label{b1}
{Let $M\subset \partial_\infty X$ be  a $H$-admissible   subset, 
$h,g \in H$   squared Type B isometries with 
  $d_*(Fix(h_*),Fix(g_*))=2$, and $W$ a wall with   $d_*(W, Fix(h_*))=1$.   
 Then    $M\cap \partial_\infty W =\phi$.}
\end{Le}
\begin{proof}
Notice  $ Fix((hgh^{-1})_*)=h_*(Fix(g_*))$     and by replacing $h$ with $h^2$ if necessary,
 we may assume  
$Fix(g_*)\cap Fix((hgh^{-1})_*)=\phi$.
Suppose $M\cap \partial_\infty W \not=\phi$. Then $M\cap \partial_\infty {h^i(W)} \not=\phi$  for every integer
 $i$ as $M$ is $H$-invariant.   By Lemma \ref{inter3}  we have 
$d_*(h^i(W), Fix(\gamma_*))\le 2$ for any integer $i$ and any  squared Type B isometry $\gamma\in H$. 
In particular,  $d_*(h^i(W), Fix(g_*))\le 2$  and 
$d_*(h^i(W), Fix((hgh^{-1})_*))\le 2$.

Since  $d_*(Fix(h_*),Fix(g_*))=2$ and    $d_*(W, Fix(h_*))=1$, Lemma \ref{atleast5}
 implies there are at most two   integers  $i$    with 
  $d_*(h^i(W), Fix(g_*))=1$.  Therefore  for all except  at most 4 integers $i$ 
 we have  
$$d_*(h^i(W), Fix(g_*))=d_*(h^i(W), Fix((hgh^{-1})_*))=2.$$   
 Fix  five distinct walls  $W_1, W_2, W_3, W_4, W_5\in \{h^i(W):i=1,2, \cdots \}$
   with 
$$d_*(W_i, Fix(g_*))=d_*(W_i, Fix((hgh^{-1})_*))=2.$$

For each $W_i$ ($i=1,\cdots, 5$),  since   $M\cap \partial_\infty {W_i} \not=\phi$
 Lemma \ref{inter2} implies there is a wall $W'_i$ with
  $$d_*(W_i, W'_i)=d_*(W'_i,Fix(g_*))=d_*(W'_i,Fix((hgh^{-1})_*))=1.$$ 
Let $\tilde W$  be the wall such  that  $h_*(\tilde W)=\tilde W$ and $d_*(W, \tilde W)=1$. 
If $W'_i=W'_j$ for some $i\not=j$,  then $\tilde WW_iW'_iW_j\tilde W$ is an injective edge
 loop of length 4, 
contradicting to Lemma \ref{atleast5}.
Suppose $W'_i\not=W'_j$ for  $i\not=j$. Since  $Fix(g_*)$ and $Fix((hgh^{-1})_*)$ are disjoint edges,
  $d_*(W'_i,Fix(g_*))=d_*(W'_i,Fix((hgh^{-1})_*))=1$   for  $i=1, \cdots,   5$
 implies there are two indices $i\not=j$, 
 and  
  walls      $W'\not=W''$ such that $g_*(W')=W'$, $(hgh^{-1})_*(W'')=W''$
   and $d_*(W'_i,W')=d_*(W'_i,W'')=d_*(W'_j,W')=d_*(W'_j,W'')=1$. Thus 
$W'_iW'W'_jW''W'_i$ is an injective edge loop of length 4 in $G_*$, a contradiction.

\end{proof}

\begin{Le}\label{b2}
{Let $M\subset \partial_\infty X$ be  a   $H$-admissible  subset,   
$h,g \in H$      squared  Type B  isometries   with   
  $d_*(Fix(h_*),Fix(g_*))=2$, and $W$ a wall with   $d_*(W, Fix(h_*))=2$.   
  Then    $M\cap \partial_\infty W =\phi$.}
\end{Le}
\begin{proof}
Since $d_*(W, Fix(h_*))=2$, Lemma \ref{inter2} implies there  is a set of  walls
 $\{W_i\}_{i=1}^N$ ($N=1$ or $2$)  with $d_*(W_i, Fix(h_*))=1$ and 
$M\cap \partial_\infty W \subset \cup_{i=1}^N (M\cap \partial_\infty {W_i})$.  
Now the lemma follows from Lemma \ref{b1}.

\end{proof}

Lemmas \ref{b0}, \ref{b1}, \ref{b2},  Proposition \ref{empadm}
 and the remark before 
Lemma \ref{a0}  
imply that  for any  two  squared  Type B  isometries $g,h\in H$, $d_*(Fix(h_*),Fix(g_*))\not=2$  holds. 
 Lemma \ref{lin} then implies 
there are   $h, h'\in H$ with
$d_*(Fix(h_*), Fix((h')_*))\ge 3$.  Now it follows from Lemma \ref{hh'} that $H$ has the TA-property. 
 The proof of Theorem \ref{sameangle}   is now complete.




 \addcontentsline{toc}{subsection}{References}

\end{document}